\documentclass{siamonline190516}

\usepackage[applemac]{inputenx}
\usepackage[english]{babel}
\usepackage[caption=false]{subfig}
\usepackage{amssymb}
\usepackage[export]{adjustbox}

\newcommand{\added}[1]{#1}
\newcommand{\removed}[1]{}

\usepackage{tikz}

\newcommand{\Colorwheel}[4]{%
        \begin{scope}%
        \clip (#1,#2) circle (#3);%
        \draw(#1,#2) node {\includegraphics[width=#4]{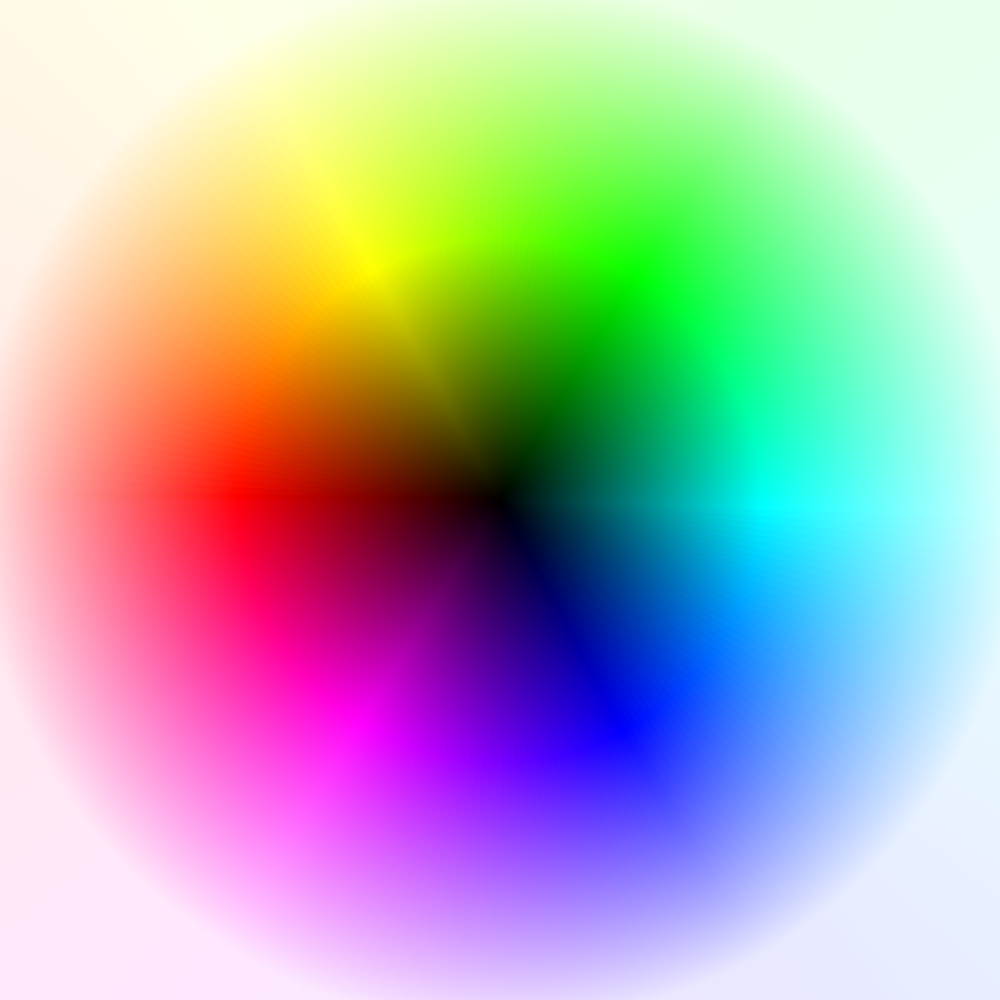}};%
        \end{scope}
        \draw(#1+#3,#2) node[anchor=west] {$0$};%
        \draw(#1-#3,#2) node[anchor=east] {$\pi$};}%
       
\newcommand{\DrawColorwheel}{%
\begin{tikzpicture}%
\Colorwheel{3cm}{2cm}{.4cm}{.8cm}%
\end{tikzpicture}}%

\usepackage{lipsum}
\usepackage{amsfonts}
\usepackage[export]{adjustbox}
\usepackage{epstopdf}
\usepackage{algorithmic}
\ifpdf
  \DeclareGraphicsExtensions{.eps,.pdf,.png,.jpg}
\else
  \DeclareGraphicsExtensions{.eps}
\fi

\numberwithin{theorem}{section}

\newcommand{\TheTitle}{Choose your path wisely: gradient descent in a Bregman distance framework}
\newcommand{\ShortTitle}{Choose your path wisely} 
\newcommand{\TheAuthors}{M. Benning, M. M. Betcke, M. J. Ehrhardt \and C.-B. Sch\"{o}nlieb.}
\newcommand{\coillength}{0.2}

\headers{\ShortTitle}{\TheAuthors}

\title{{\TheTitle}\thanks{Submitted to the editors DATE.
\funding{This work was funded by the Leverhulme Trust Early Career Fellowship 'Learning from mistakes: a supervised feedback-loop for imaging applications', the Isaac Newton Trust, the Engineering and Physical Sciences Research Council (EPSRC) 'EP/K009745/1', the Leverhulme Trust project 'Breaking the non-convexity barrier', the EPSRC grant 'EP/M00483X/1', the EPSRC centre 'EP/N014588/1', the Cantab Capital Institute for the Mathematics of Information and CHiPS (Horizon 2020 RISE project grant).}}}

\author{
  Martin Benning\thanks{School of Mathematical Sciences, Queen Mary University of London, UK (\email{m.benning@qmul.ac.uk}).}
    \and
  Marta M. Betcke\thanks{Department of Computer Science, University College London, UK (\email{m.betcke@ucl.ac.uk}).}    
    \and
   Matthias J. Ehrhardt\thanks{Institute for Mathematical Innovation, University of Bath, UK (\email{m.ehrhardt@bath.ac.uk}).}
  \and
  Carola-Bibiane Sch\"{o}nlieb\thanks{Department of Applied Mathematics and Theoretical Physics, University of Cambridge, UK (\email{cbs31@cam.ac.uk}).}
}

\usepackage{amsopn}

\DeclareMathOperator*{\argmin}{\arg \min}
\DeclareMathOperator{\dom}{dom}
\DeclareMathOperator{\crit}{crit}

\DeclareMathOperator{\dist}{dist}

\newcommand{\FunConvex}{\Gamma_0}
\newcommand{\FunSmooth}[1]{\mathcal S_{#1}}

\newcommand{\FunProblem}[1][L]{\Psi_{#1}}

\newcommand{\C}{\mathbb{C}}
\newcommand{\R}{\mathbb{R}}
\newcommand{\N}{\mathbb{N}}
\DeclareMathOperator{\uc}{\R^n}

\newcommand{\taumin}{\tau^{\text{min}}}

\newtheorem{remark}{Remark}

\newcommand{\changed}[1]{{\color{black}#1}}

\begin{document}

\maketitle

\begin{abstract}
We propose an extension of a special form of gradient descent --- in the literature known as linearised Bregman iteration --- to a larger class of non-convex functions. We replace the classical (squared) two norm metric in the gradient descent setting with a generalised Bregman distance, based on a proper, convex and lower semi-continuous function. The algorithm's global convergence is proven for functions that satisfy the Kurdyka-\L ojasiewicz property. Examples illustrate that features of different scale are being introduced throughout the iteration, transitioning from coarse to fine. This coarse-to-fine approach with respect to scale allows to recover solutions of non-convex optimisation problems that are superior to those obtained with conventional gradient descent, or even projected and proximal gradient descent. The effectiveness of the linearised Bregman iteration in combination with early stopping is illustrated for the applications of parallel magnetic resonance imaging, blind deconvolution as well as image classification with neural networks.
\end{abstract}



\section{Introduction}
Non-convex optimisation methods are indispensable mathematical tools for a large variety of applications \cite{nocedal2006numerical}. For differentiable objectives, first-order methods such as gradient descent have proven to be useful tools in all kinds of scenarios. Throughout the last decade, however, there has been an increasing interest in first-order methods for non-convex and non-smooth objectives. These methods range from forward-backward, respectively proximal-type, schemes \cite{attouch2009convergence,attouch2010proximal,attouch2013convergence,Bonettini2015,Bonettini2016}, over linearised proximal schemes \cite{xu2013block,bolte2014proximal,xu2017globally,nikolova2017alternating}, to inertial methods \cite{ochs2014ipiano,pock2016inertial}, primal-dual algorithms \cite{valkonen2014primal,li2015global,moeller2015variational,benning2015preconditioned}, scaled gradient projection methods \cite{Prato.ea2016} and non-smooth Gau\ss-Newton extensions \cite{drusvyatskiy2016nonsmooth,ochs2017non}.

In this paper, we follow a different approach of incorporating non-smoothness into first-order methods for non-convex problems. We present a direct generalisation of gradient descent, first introduced in \cite{benning2016gradient}, where the usual squared two-norm metric that penalises the gap of two subsequent iterates is being replaced by a potentially non-smooth distance term. This distance term is given in form of a generalised Bregman distance \cite{bregman1967relaxation,burger2013adaptive,osher2016sparse}, where the underlying function is proper, lower semi-continuous and convex, but not necessarily smooth. If the underlying function is a Legendre function (see \cite[Section 26]{rockafellar1970convex} and \cite{bauschke2001essential}), the proposed generalisation basically coincides with the recently proposed non-convex extension of the Bregman proximal gradient method \cite{bolte2017first}. In the more general case, the proposed method is a generalisation of the so-called linearised Bregman iteration \cite{darbon2007fast,yin2008bregman,cai2009linearized,cai2009convergence} to non-convex data fidelities.

Motivated by inverse scale space methods (cf. \cite{burger2006nonlinear,burger2013adaptive,osher2016sparse}), the use of non-smooth Bregman distances for the penalisation of the iterates gap allows to control the scale of features present in the individual iterates. Replacing the squared two-norm, for instance, with a squared two-norm plus the Bregman distance w.r.t. a one-norm leads to very sparse initial iterates, with iterates becoming more dense throughout the course of the iteration. This control of scale, i.e. the slow evolution from iterates with coarse structures to iterates with fine structures, can help tp overcome unwanted minima of a non-convex objective, as we are going to demonstrate with an example in Section \ref{sec:motivation}. This is in stark contrast to many of the non-smooth, non-convex first-order approaches mentioned above, where the methods are often initialised with random inputs that become more regular throughout the iteration.

Our main contributions of this paper are the generalisation of the linearised Bregman iteration to non-convex functions, a detailed convergence analysis of the proposed method as well as the presentation of numerical results that demonstrate that the use of coarse-to-fine scale space approaches in the context of non-convex optimisation can lead to superior solutions. 

The outline of the paper is as follows. Based on the non-convex problem of blind deconvolution, we first give a motivation in Section \ref{sec:motivation} of why a coarse-to-fine approach in terms of scale can indeed lead to superior solutions of non-convex optimisation problems. We then recall key concepts of convex and non-convex analysis that are needed throughout the paper in Section \ref{sec:mathprelim}. Subsequently, we define the extension of the linearised Bregman iteration for non-convex functions in Section \ref{sec:linbreg}. Then, motivated by the informal convergence recipe of Bolte et al. \cite[Section 3.2]{bolte2014proximal} we show a global convergence result in Section \ref{sec:globconv}, which concludes the theoretical part. We conclude with the modelling of the applications of parallel Magnetic Resonance Imaging (MRI), blind deconvolution and image classification in Section \ref{sec:apps}, followed by corresponding numerical results in Section \ref{sec:numresults} as well as conclusions and outlook in Section \ref{sec:concoutlook}.

\section{Motivation}\label{sec:motivation}
\begin{figure}[!t]\label{fig:pixelblinddeconvolution}
\begin{center}
\subfloat[Original image $\hat{u}$]{\includegraphics[width=0.32\textwidth]{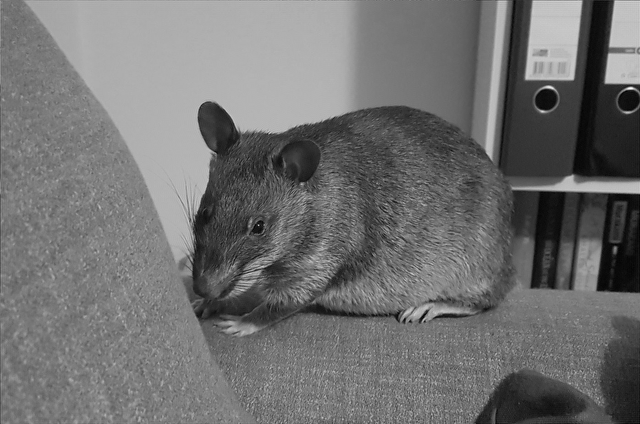}\label{subfig:introduction1}}\vspace{0.05cm}
\subfloat[$f$ and $\hat{h}$]{\includegraphics[width=0.32\textwidth]{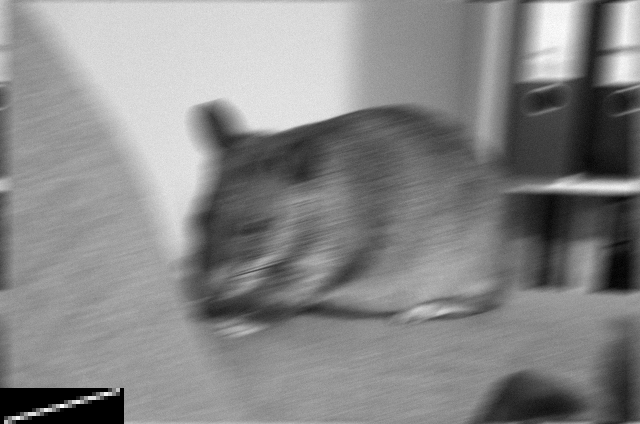}\label{subfig:introduction2}}\vspace{0.05cm}
\subfloat[Projected gradient descent]{\includegraphics[width=0.32\textwidth]{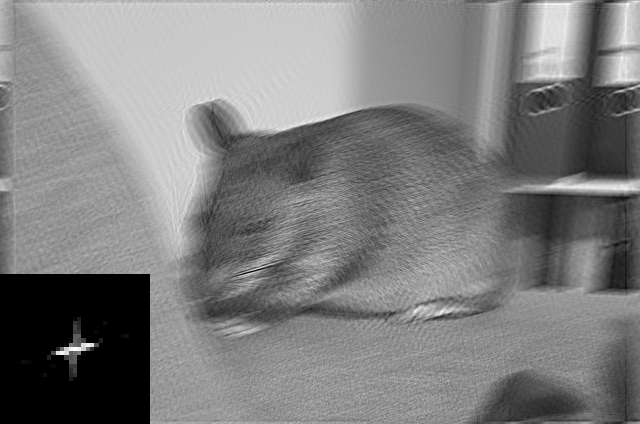}\label{subfig:introduction3}}\vspace{0.05cm}\\
\subfloat[$\alpha = 10^{-3}$]{\includegraphics[width=0.32\textwidth]{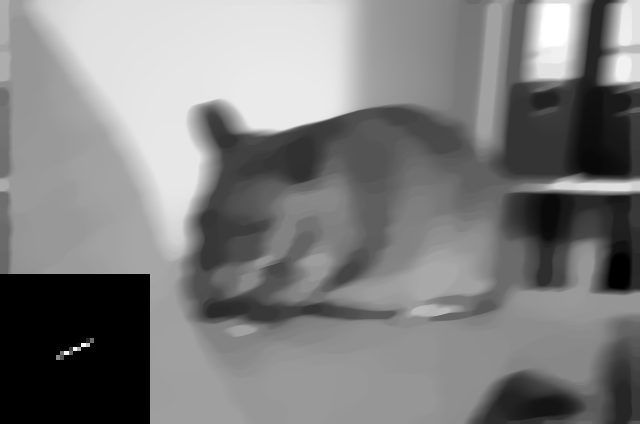}\label{subfig:introduction5}}\vspace{0.05cm}
\subfloat[$\alpha = 10^{-4}$]{\includegraphics[width=0.32\textwidth]{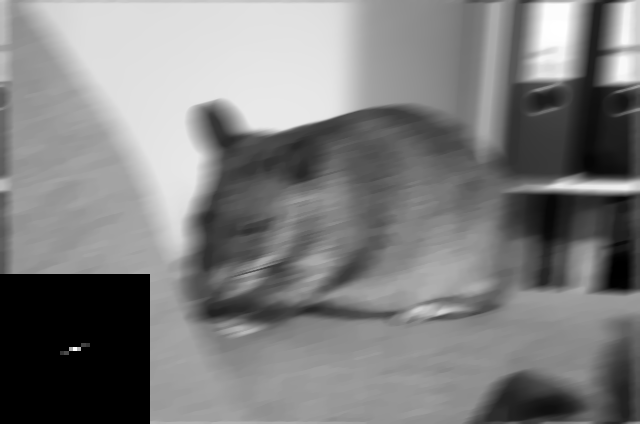}\label{subfig:introduction6}}\vspace{0.05cm}
\end{center}
\caption{\textbf{Standard approaches for blind deconvolution.} Figure \ref{subfig:introduction1} shows the image $\hat{u}$ of Pixel the Gambian pouched rat, courtesy of Monique Boddington. Figure \ref{subfig:introduction2} shows a motion-blurred version $f$ of that same image; the corresponding convolution kernel $\hat{h}$ is depicted in the bottom left corner. Figure \ref{subfig:introduction3} visualises the reconstruction of the image and the convolution kernel obtained with the projected gradient descent method \eqref{subeq:pgd}. In Figure \ref{subfig:introduction5} we see the result of gradient descent method \eqref{eq:proxgd} for $\alpha = 10^{-3}$, whereas Figure \ref{subfig:introduction6} shows the result of \eqref{eq:proxgd} for the choice $\alpha = 10^{-4}$.}
\end{figure}

We want to motivate the use of the linearised Bregman iteration for non-convex optimisation problems with the example of blind deconvolution. In blind (image) deconvolution the goal is to recover an unknown image $u$ from a blurred and usually noisy image $f$. Assuming that the degradation is the same for each pixel, the problem of blind deconvolution can be modelled as the minimisation of the energy
\begin{align}
E_{1}(u, h) := \underbrace{\frac{1}{2} \| u \ast h - f \|_2^2}_{=: F(u, h)} + \chi_{C}(h) \, ,\label{eq:blinddeconv}
\end{align}
with respect to the arguments $u \in \R^n$ and $h \in \R^r$. Here $\ast$ denotes a discrete convolution operator, and $\chi_{C},$
is the characteristic function
\begin{align*}
\chi_{C}(h) := \begin{cases} 0 & h \in C\\ \infty & h \not\in C \end{cases} \, ,
\end{align*}
defined over the simplex constraint set
\begin{align*}
C := \left\{ h \in \mathbb{R}^r \ \left|  \  \sum_{j = 1}^r h_j = 1, \ h_j \geq 0, \ \forall j \in \{1, \ldots, r \} \right. \right\} \, .
\end{align*}
Even with data $f$ in the range of the non-linear convolution operator, i.e. $f = \hat{u} \ast \hat{h}$ for some $\hat{u} \in \R^n$ with $\hat{h} \in C$, it is usually still fairly challenging to recover $\hat{u}$ and $\hat{h}$ as solutions of \eqref{eq:blinddeconv}. A possible reason for this could be that \eqref{eq:blinddeconv} is an invex function on $\R^n \times C$, where every stationary point is already a global minimum. If we simply try to recover $\hat{u}$ and $\hat{h}$ via projected gradient descent, we usually require an initial point in the neighbourhood of $(\hat{u}, \hat{h})$ in order to converge to that point. We want to illustrate this with a concrete example. Assume we are given an image $\hat{u}$ and a convolution kernel $\hat{h}$ as depicted in Figure \ref{fig:pixelblinddeconvolution}, and $f = \hat{u} \ast \hat{h}$ is as shown in Figure \ref{subfig:introduction2}.
Minimising \eqref{eq:blinddeconv} via projected gradient descent leads to the following procedure:%
\begin{subequations}%
\begin{align}%
u^{k + 1} &= u^k - \tau^k \, \partial_u \, F(u^k, h^k) \, ,\label{subeq:pgd1}\\%
h^{k + 1} &= \text{proj}_{C}\left(h^k - \tau^k \, \partial_h \, F(u^k, h^k) \right) \, ,%
\end{align}\label{subeq:pgd}%
\end{subequations}%
\changed{where $\text{proj}_{C}$ denotes the projection onto the convex set $C$}. If we initialise with $u^0 = (0, \ldots, 0)^T$ and $h^0 = (1, \ldots, 1)^T/r$, set $\tau^0 = 1$, update $\tau^k$ via backtracking to ensure a monotonic decrease of the energy $E_1$, and iterate \eqref{subeq:pgd} for 3500 iterations, we obtain the reconstructions visualised in Figure \ref{subfig:introduction3}. Even without any noise present in the data $f$, the algorithm converges to a solution very different from $\hat{u}$ and $\hat{h}$. This is not necessarily surprising as we do not impose any regularity on the image. We can try to overcome this issue by modifying \eqref{eq:blinddeconv} as follows:
\begin{align}
E_2(u, h) &:= F(u, h) + \chi_{C}(h) + \alpha \text{TV}(u) \, , \label{eq:blinddeconvmod}\\
&= E_1(u, h) + \alpha \text{TV}(u) \, . \nonumber
\end{align}
Here TV denotes the discretised total variation, i.e.%
\begin{align*}%
\text{TV}(u) := \| | \nabla u | \|_1 \, ,%
\end{align*}%
where $\nabla:\R^n \rightarrow \R^{2n}$ is a (forward) finite difference discretisation of the gradient operator, $| \cdot |$ the Euclidean vector norm and $\| \cdot \|_1$ the one-norm, and $\alpha$ is a positive scalar. The minimisation of \eqref{eq:blinddeconvmod} can easily be carried out by the proximal gradient descent method, also known as forward-backward splitting \cite{lions1979splitting}, which is a minor modification of the projected gradient method \cite{goldstein1964convex,goldstein1967constructive,bertsekas1976goldstein} to more general proximal mappings. In the context of minimising \eqref{eq:blinddeconvmod}, the proximal gradient method reads as%
\begin{subequations}%
\begin{align}%
u^{k + 1} &= (I + \alpha \partial \text{TV})^{-1}(u^k - \tau^k \, \partial_u \, F(u^k, h^k) ) \, ,\\%
h^{k + 1} &= \text{proj}_{C}\left(h^k - \tau^k \, \partial_h \, F(u^k, h^k) \right)  \, ,%
\end{align}\label{eq:proxgd}%
\end{subequations}%
where $(I + \alpha \partial \text{TV})^{-1}$ denotes the proximal mapping \cite{moreau1962decomposition,moreau1965proximite} with respect to the total variation, i.e.
\begin{align}
(I + \alpha \partial \text{TV})^{-1}(z) := \argmin_{u \in \R^n} \left\{ \frac{1}{2} \| u - z \|_2^2 + \alpha \text{TV}(u) \right\} \, .\label{eq:tvprox}
\end{align}
It is straight-forward to solve \eqref{eq:tvprox} for a given argument with numerical methods such as the (accelerated) primal-dual hybrid gradient method (cf. \cite{zhu2008efficient,pock2009algorithm,esser2010general,chambolle2011first,chambolle2016introduction}) up to sufficient numerical accuracy. If we then evaluate 3000 iterations of \eqref{eq:proxgd} for $\alpha \in \{ 10^{-3}, 10^{-4} \}$ with the same initial values that we used for the projected gradient method, we obtain the results visualised in Figure \ref{fig:pixelblinddeconvolution}. 
\begin{figure}[!t]\label{fig:introductioniterates}
\begin{center}
\subfloat[1st iterate]{\includegraphics[width=0.32\textwidth]{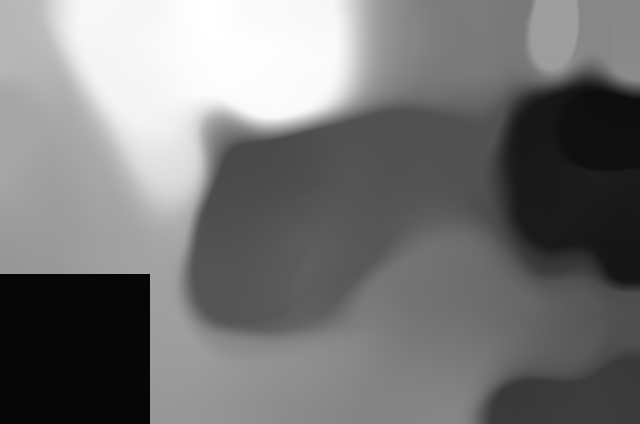}\label{subfig:introduction7}}\vspace{0.05cm}
\subfloat[10th iterate]{\includegraphics[width=0.32\textwidth]{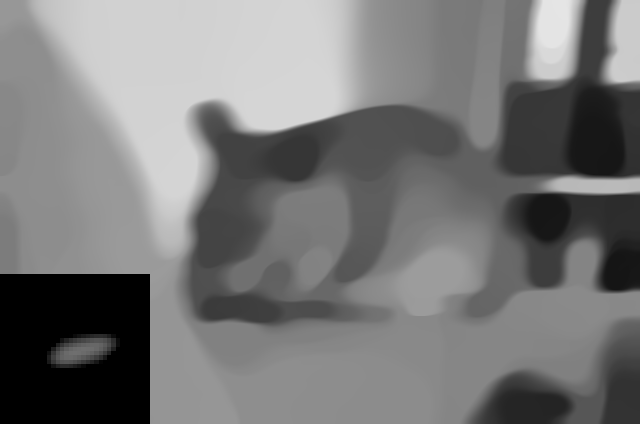}\label{subfig:introduction8}}\vspace{0.05cm}
\subfloat[50th iterate]{\includegraphics[width=0.32\textwidth]{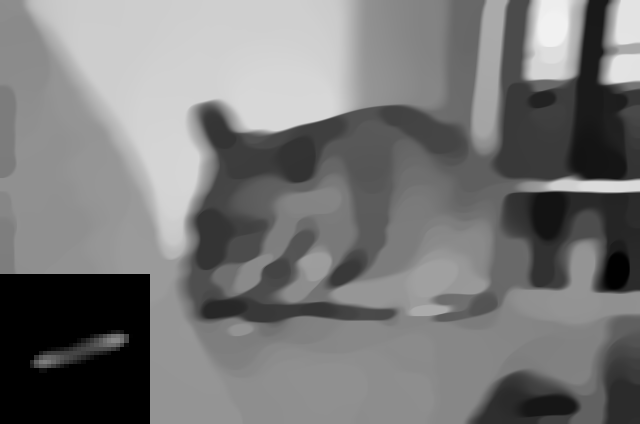}\label{subfig:introduction9}}\\
\subfloat[500th iterate]{\includegraphics[width=0.32\textwidth]{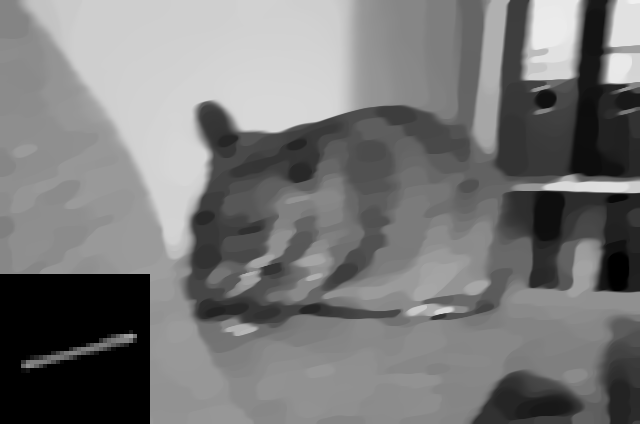}\label{subfig:introduction10}}\vspace{0.05cm}
\subfloat[1500th iterate]{\includegraphics[width=0.32\textwidth]{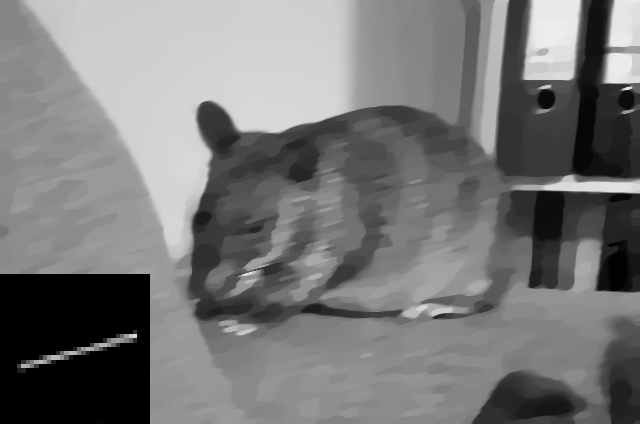}\label{subfig:introduction11}}\vspace{0.05cm}
\subfloat[3000th iterate]{\includegraphics[width=0.32\textwidth]{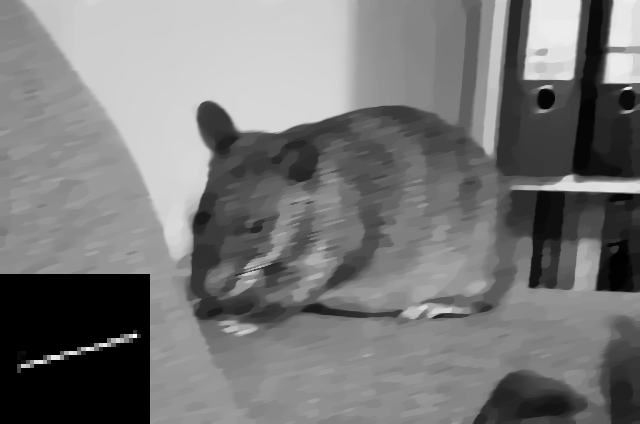}\label{subfig:introduction12}}
\end{center}
\caption{\textbf{Proposed approach for blind deconvolution.} Figure \ref{fig:introductioniterates} shows several iterates of the linearised Bregman iteration \eqref{eq:introductionlinbregman} for the choice $\alpha = 0.05$. The strong initial effect of the total variation regularisation enables the algorithm to converge to a solution close to $\hat{u}$ and $\hat{h}$.} 
\end{figure}
We observe that for the larger choice of $\alpha = 10^{-3}$ we obtain a better reconstruction of the convolution kernel, but at the cost of a reconstructed image that is very cartoon-like. Reducing the parameter $\alpha$ to $\alpha = 10^{-4}$ reduces the impact of the total variation regularisation; however, the reconstructed image then remains fairly blurry and the reconstructed convolution kernel is closer to a Dirac delta.

The reason for this is that the total variation-based model \eqref{eq:blinddeconvmod} is basically not suitable for deconvolution tasks. Blurred images generally have a smaller total variation compared to their sharp counterparts, hence it is easier to minimise the energy in \eqref{eq:blinddeconvmod} by recovering a kernel close to a Dirac delta and a smoothed version of the blurry image in order to reduce the total variation. 

We therefore want to use an alternative approach that is different to the two approaches presented above. We do observe from the proximal gradient example that a larger regularisation parameter seems to work better for a more accurate reconstruction of the convolution kernel (at the cost of a rather cartoon-like image). The explanation for this is that image features at a relatively coarse scale have to be adjusted to minimise the data fit, forcing the convolution kernel to correct for this. It therefore seems reasonable to find a minimiser of \eqref{eq:blinddeconv} with a scale-space approach, changing from coarse to fine scales over the course of the iteration. Specifically, we propose to use a variant of the linearised Bregman iteration adopted to minimising non-convex problems such as the minimisation of the function $E_1$ as defined in \eqref{eq:blinddeconv}. For the choice of $E_1$ in \eqref{eq:blinddeconv}, this method reads as
\begin{subequations}%
\begin{align}%
u^{k + 1} &= \argmin_{u \in \R^n} \left\{ \frac{1}{2} \| u - u^k \|^2 + \tau^k \left( \alpha D_{\text{TV}}^{q^k}(u, u^k) + \langle \partial_u \, F(u^k, h^k), u \rangle \right) \right\} \, ,\label{subeq:introduction2}\\%
q^{k + 1} &=  q^k - \frac{1}{\tau^k \alpha} \left( u^k - u^{k + 1} - \tau^k \, \partial_u \, F(u^k, h^k) \right) \, ,\\%
h^{k + 1} &= \text{proj}_C\left(h^k - \tau^k \, \partial_h \, F(u^k, h^k) \right) \, .%
\end{align}\label{eq:introductionlinbregman}%
\end{subequations}%
Here $q^k \in \partial \text{TV}(u^k)$ denotes a subgradient of $\text{TV}$ at $u^k$, $\alpha \geq 0$ is a scalar and $D^{q^k}_{\text{TV}}(u^{k + 1}, u^k)$ is the generalised Bregman distance \cite{bregman1967relaxation} with respect to the total variation, i.e.
\begin{align*}
D^{q^k}_{\text{TV}}(u^{k + 1}, u^k) = \text{TV}(u^{k + 1}) - \text{TV}(u^k) - \langle q^k, u^{k + 1} - u^k \rangle \, ,
\end{align*}
for a subgradient $q^k \in \partial \text{TV}(u^k)$. Note that \eqref{eq:introductionlinbregman} reduces to the projected gradient method \eqref{subeq:pgd} for the choice $\alpha = 0$.

Replacing the total variation semi-norm in \eqref{eq:proxgd} with its Bregman distance yields an iterative scale-space method that changes the influence of the total variation regularisation throughout the course of the iteration. With a larger parameter $\alpha$, the initial iterates have a very low total variation and contain only coarse features. Throughout the iteration, features of finer and finer scale are introduced. We have visualised several iterates of \eqref{eq:introductionlinbregman} for the choice $\alpha = 0.05$ in Figure \ref{fig:introductioniterates} to demonstrate this phenomenon.

We observe that this modification of projected gradient descent enables us to converge to minimisers of $E_1$ as defined in \eqref{eq:blinddeconv} that are fairly close to the original choices of $\hat{u}$ and $\hat{h}$. Hence, the choice of Bregman distance strongly affects the outcome of the iteration procedure and can be used to guide the iterates towards more desirable outcomes. 

Obviously real data is never in the range of the forward model, and in that case we do not want to converge to a minimiser of $E_1$. However, we can still apply the linearised Bregman iteration in combination with early stopping in order to produce superior results compared to projected or proximal gradient descent, which we will further demonstrate in Section \ref{sec:apps} and Section \ref{sec:numresults}. Prior to this, we provide a comprehensive convergence analysis of the linearised Bregman iteration in the Sections \ref{sec:linbreg} and \ref{sec:globconv}.

\section{Linearised Bregman iteration for non-convex problems}\label{sec:linbreg}
We are interested in the minimisation of functions $E \in \FunSmooth{L}$, where $\FunSmooth{L}$ is defined in Definition \ref{def:smoothness} in the appendix. We want to emphasise that the function $E$ does not necessarily have to be convex. In order for the minimisation of $E$ to make sense, we have to introduce some additional assumptions for this function first. From now on we assume $E \in \FunProblem$, with $\FunProblem$ being defined as
\begin{align*}
\FunProblem := \left\{ E \in \FunSmooth{L} \left| \begin{array}{c} \text{$E$ has bounded level sets} \\ \text{$E$ is bounded from below} \end{array} \right. \right\} \, .
\end{align*}


\noindent We further recall the definition of the set of critical points of $E$, i.e.
\begin{align}
\crit(E) := \left\{ u \in \dom(E) \, \left| \, \nabla E(u) = 0 \right. \right\} \, .\label{eq:crit}
\end{align}
The requirements on $E$ ensure that sequences $\{ u^k \}_{k \in \N}$ are already bounded if the sequences $\{ E(u^k) \}_{k \in \N}$ are bounded, that an infimum exists and that the set of critical points is non-empty.

\begin{algorithm}[t]
\caption{Generalised linearised Bregman iteration for minimising $E$}
\label{alg:linbreg}
\begin{algorithmic}
\STATE{Initialise $\{ \tau^k \}_{k \in \N}$, $u^0$ and $p^0 \in \partial J(u^0)$}
\FOR{$k = 0, 1, \ldots$}
\STATE{Compute $u^{k + 1} = \argmin_{u \in \uc} \left\{ \tau^k \langle u - u^k, \nabla E(u^k) \rangle + D_J^{p^k}(u, u^k) \right\}$}
\STATE{Compute $p^{k +1} = p^k - \tau^k \nabla E(u^k)$}
\ENDFOR
\end{algorithmic}
\end{algorithm}

We want to minimise $E$ iteratively in a way that allows us to follow solution paths of different regularity. This regularity will be induced by an additional function $J \in \FunConvex$, where $\FunConvex$ is defined in the appendix. Precisely, we approach the minimisation of $E$ via the linearised Bregman iteration
\begin{subequations}\label{eq:linbregman}
\begin{align}
u^{k + 1} &= \argmin_{u \in \uc} \left\{ \tau^k \langle \nabla E(u^k), u - u^k\rangle + D_J^{p^k}(u, u^k) \right\} \, ,\label{eq:linbregman1} \\
p^{k +1} &= p^k - \tau^k \nabla E(u^k) \, ,\label{eq:linbregman2}
\end{align}
\end{subequations}
for $k \in \N$, a sequence of positive parameters $\{ \tau^k \}_{k \in \N}$ and initial values $u^0$ and $p^0$ with $p^0 \in \partial J(u^0)$. Here $\partial J$ denotes the subdifferential; we refer to the appendix for its definition. Note that \eqref{eq:linbregman2} is simply the optimality condition of \eqref{eq:linbregman1}. If $J$ is differentiable, $\partial J$ is single-valued and we do not have to compute \eqref{eq:linbregman2} as we do not need to pick a specific element from the set. However, if $\partial J$ is multivalued, \eqref{eq:linbregman} guarantees $p^{k+1} \in \partial J(u^{k+1})$ for all $k \in \N$. This general form of linearised Bregman iteration for the minimisation of non-convex functions is summed up in Algorithm \ref{alg:linbreg}.
\begin{remark}
For $J(u) = \frac{1}{2} \| u \|^2$, \eqref{eq:linbregman} (and therefore also Algorithm \ref{alg:linbreg}) reduces to classical gradient descent. Hence, the linearised Bregman iteration is indeed a generalisation of gradient descent.
\end{remark}

Based on what has become known as the Bregman iteration \cite{censor1992proximal,teboulle1992entropic,eckstein1993nonlinear,kiwiel1997proximal,osher2005iterative}, the linearised Bregman iteration has initially been proposed in \cite{darbon2007fast} for the computation of sparse solutions of underdetermined linear systems of equations. It has been extensively studied in this context (cf. \cite{yin2008bregman,cai2009linearized,cai2009convergence}) and also in the context of the minimisation of more general convex functions (see \cite{yin2010analysis}). It is also closely linked to (linearised variants of) the alternating direction method of multipliers (ADMM) \cite{gabay1983chapter}, as well as generalisations to non-quadratic Bregman distances \cite{wang2014bregman}. It has further been analysed in the context of non-linear inverse problems in \cite{Bachmayr2009}. In \cite{benning2016gradient}, the linearised Bregman iteration has been studied in the context of minimising general smooth but non-convex functions. Algorithm \ref{alg:linbreg} allows us to control the scale of the iterates, depending on the choice of $J$. Note that we can also reformulate \eqref{eq:linbregman1} as follows:
\begin{align}
u^{k + 1} &= \argmin_{u \in  \uc} \left\{  \tau^k \left\langle \nabla E(u^k) - \frac{1}{\tau^k} p^k, u - u^k \right\rangle + J(u)  \right\} \, \text{.} \label{eq:linbregman1:reformulated}
\end{align}
In order to ensure that a solution of Update \eqref{eq:linbregman1:reformulated} (respectively \eqref{eq:linbregman1}) exists, we choose $J$ such that $J(u) + \tau^k \langle u^\ast , u\rangle$ is coercive 
for all $u^\ast \in \uc$. In particular, we choose $J$ to be of the form $J_k := \frac{1}{2} \| \cdot \|^2 + \tau^k R$, where $R \in \FunConvex$. For this choice the iterates \eqref{eq:linbregman} read as
\begin{subequations}\label{eq:linbregman:special}
\begin{align}
u^{k + 1} {} = {} &\argmin_{u \in \uc} \left\{ \tau^k \left( \langle \nabla E(u^k), u - u^k \rangle + D_R^{q^k}(u, u^k) \right) + \frac{1}{2} \| u - u^k \|^2 \right\} \, , \nonumber\\
{} = {} &\left(I + \tau^k \partial R \right)^{-1}\left( u^k + \tau^k \left( q^k - \nabla E(u^k) \right) \right)\, ,\label{eq:linbregman:special:primal}\\
q^{k + 1} {} = {} &q^k - \frac{1}{\tau^k} \left( u^{k + 1} - u^k + \tau^k \nabla E(u^k) \right) \, , \label{eq:linbregman:special:dual}
\end{align}
\end{subequations}
for $q^k \in \partial R(u^k)$. Note that \eqref{eq:linbregman:special:dual} can be written as
\begin{align}
q^{k + 1} = q^0 - \sum_{n = 0}^k \left[ \frac{1}{\tau^n} ( u^{n + 1} - u^n ) \right] - \sum_{n = 0}^k \nabla E(u^n) \, ,\label{eq:sumdualupdate}
\end{align}
and hence, for constant stepsize $\tau^k = \tau$ \eqref{eq:linbregman:special:primal} simplifies to
\begin{align}
\begin{split}
u^{k + 1} = \left(I + \tau \partial R \right)^{-1}\left( u^0 + \tau q^0 - \tau \sum_{n = 0}^k \nabla E(u^n) \right)
\end{split} \, .\label{eq:sumprimalupdate}
\end{align}
Equations \eqref{eq:linbregman:special} are summarised in Algorithm \ref{alg:linbregspecial}. Note that both Algorithm \ref{alg:linbregspecial} and Equation \eqref{eq:sumprimalupdate} demonstrate that this specialised linearised Bregman iteration is indeed different to proximal gradient descent, for which one iterate reads $u^{k + 1} = \left(I + \tau \partial R \right)^{-1}\left( u^k - \tau \nabla E(u^k) \right)$. \added{Instead, from Equation \eqref{eq:linbregman:special:primal} we observe that one computes a subgradient descent step in the direction of the subgradient of $E - R$, followed by an application of the proximal step with respect to $R$.}

In the following we prove decrease properties and a global convergence result for Algorithm \ref{alg:linbregspecial}.

\begin{algorithm}[t]
\caption{Specialised linearised Bregman iteration for minimising $E$}
\label{alg:linbregspecial}
\begin{algorithmic}
\STATE{Initialise $\{ \tau^k \}_{k \in \N}$, $u^0$ and $q^0 \in \partial R(u^0)$}
\FOR{$k = 0, 1, \ldots$}
\STATE{Get $u^{k + 1} = \left(I + \tau^k \partial R \right)^{-1}\left( u^k + \tau^k \left( q^k - \nabla E(u^k) \right) \right)$}
\STATE{Compute $q^{k + 1} = q^k - \frac{1}{\tau^k} \left( u^{k + 1} - u^k + \tau^k \nabla E(u^k) \right)$}
\ENDFOR
\end{algorithmic}
\end{algorithm}

\section{A global convergence result for Algorithm \ref{alg:linbregspecial}}\label{sec:globconv}
The convergence analysis is inspired by the global convergence recipe of \cite{bolte2014proximal}. It is an extension to a class of non-smooth surrogate functions for which a tailored convergence analysis is presented that utilises the convexity of $R$. We begin our analysis of Algorithm \ref{alg:linbregspecial} by showing a sufficient decrease property of the surrogate function and a subgradient bound by the (primal) iterates gap. In order to do so, we first define the following surrogate function for $E$.

\begin{definition}[Surrogate objective]\label{def:surrogate}
Assume $E \in \FunProblem$ and $R \in \FunConvex$. Then we define a surrogate function $F:\R^n \times \R^n \rightarrow \R \cup \{ \infty \}$ as
\begin{align}
F(x, y) := E(x) + R(x) + R^\ast(y) - \langle x, y \rangle .\label{eq:surrogate}
\end{align}
Here $R^\ast$ denotes the convex conjugate of $R$ as defined in Definition \ref{def:convconj} in the appendix. 
\end{definition}
Note that based on Remark \ref{rem:bregdis} in the appendix, the surrogate function \eqref{eq:surrogate} satisfies
\begin{align*}
F(x, y) := E(x) + D_R^{y}(x, z) \, ,
\end{align*}
for any $z \in \partial R^\ast(y)$, which implies $F(x, y) \geq E(x)$ for all $x, y \in \uc$. Before we continue, we want to introduce the concise notation $s^k := (u^k, q^{k - 1})$ for all $k \in \N$, such that $F(s^k) = F(u^k, q^{k-1})$. 
With the following lemma we prove a sufficient decrease property of the surrogate energy \eqref{eq:surrogate} for subsequent iterates.
\begin{lemma}[Sufficient decrease property]\label{thm:sufficientdecrease}
Assume $E \in \FunProblem$ and $R \in \FunConvex$. Further, suppose that the stepsize $\tau^k$ satisfies the condition
\begin{align}
0 < \tau^k \leq \frac{2}{L + 2 \rho_1} \, ,\label{eq:stepsizebound}
\end{align}
for some $\rho_1 > 0$ and all $k \in \N$. Then the iterates of Algorithm \ref{alg:linbregspecial}
 satisfy the descent estimate
\begin{align}
F(s^{k + 1}) + \rho_1 \| u^{k + 1} -  u^k \|^2  \leq F(s^k) \, ,\label{eq:suffdecrease}
\end{align}
for $s^k := (u^k, q^{k - 1})$ and $F$ as defined in \eqref{eq:surrogate}. In addition, we observe
\begin{align}
\lim_{k \rightarrow \infty} \| u^{k + 1} -  u^k \|^2 = 0 \quad \text{as well as} \quad \lim_{k \rightarrow \infty} D_R^{\text{symm}}(u^{k + 1}, u^k) = 0 \, . \label{eq:lim}
\end{align}
\begin{proof}
First of all, we compute
\begin{align*}
\tau^k \left( \nabla E(u^k) + q^{k + 1} - q^k \right) + u^{k + 1} - u^k = 0
\end{align*}
as the optimality condition of \eqref{eq:linbregman:special:primal}, which is also the rearranged update formula \eqref{eq:linbregman:special:dual} as mentioned earlier (for $q^{k + 1} \in \partial R(u^{k + 1})$). Taking the inner product with $u^{k + 1} - u^k$ therefore yields
\begin{align}
-\langle \nabla E(u^k), u^{k + 1} - u^k \rangle = \frac{1}{\tau^k} \| u^{k + 1} -  u^k \|^2 +  D^{\text{symm}}_R(u^{k + 1}, u^k) \, .\label{eq:intermedeq}
\end{align}
Due to the Lipschitz-continuity of the gradient of $E$ we can use \eqref{eq:lipschitzest3} from the appendix and further estimate
\begin{align*}
E(u^{k + 1}) \leq E(u^k) + \langle \nabla E(u^k), u^{k + 1} - u^k\rangle + \frac{L}{2} \| u^{k + 1} - u^k \|^2 \, .
\end{align*}
Together with \eqref{eq:intermedeq} and the stepsize bound \eqref{eq:stepsizebound} we therefore obtain the estimate
\begin{align}
E(u^{k + 1}) +  D_R^{\text{symm}}(u^{k+1}, u^k) + \rho_1 \| u^{k + 1} -  u^k \|^2 \leq E(u^k) \, . \label{eq:desc:E}
\end{align}
Adding $ D_R^{q^{k - 1}}(u^k, u^{k - 1})$ to both sides of the inequality then allows us to conclude
\begin{align*}
&F(s^{k + 1}) +  D_R^{q^{k + 1}}(u^k, u^{k + 1}) + D_R^{q^{k - 1}}(u^k, u^{k - 1}) + \rho_1 \| u^{k + 1} -  u^k \|^2 \\
{} \leq {}  &F(s^{k}) \, .
\end{align*}
Due to the non-negativity of $D_R^{q^{k + 1}}(u^k, u^{k + 1})$ and $D_R^{q^{k - 1}}(u^k, u^{k - 1})$, we have verified \eqref{eq:suffdecrease}. Moreover, summing up \eqref{eq:desc:E} over $k = 0, \dots, N$ yields
\begin{align*}
\sum_{k = 0}^N \left[ \rho_1  \| u^{k + 1} -  u^k \|^2 +  D_R^\text{symm}(u^{k + 1}, u^k) \right] &\leq \sum_{k = 0}^N E(u^k) - E(u^{k + 1}) \, \text{,}\\
&= E(u^0) - E(u^{N + 1}) \, \text{,}\\
&\leq E(u^0) - \inf_u E(u) < \infty \, \text{.}
\end{align*}
Taking the limit $N \rightarrow \infty$ therefore implies
\begin{align*}
\sum_{k = 0}^\infty \left[ \rho_1 \| u^{k + 1} -  u^k \|^2 + D_R^\text{symm}(u^{k + 1}, u^k) \right] < \infty \, \text{,}
\end{align*}
and thus \eqref{eq:lim}, due to $\rho_1 > 0$.
\end{proof}
\end{lemma}


\begin{remark}\label{rem:boundedsequence}
As Lemma \ref{thm:sufficientdecrease} implies the monotonic decrease $F(s^{k + 1}) \leq F(s^k)$, we already know that the sequence $\{ F(s^k) \}_{k \in \N}$ is bounded from above. It is also bounded from below, since $F(s^k) \geq E(u^k) \geq \inf_u E(u) > - \infty$, due to $E \in \FunProblem$.
\end{remark}
It is worth mentioning that the name sufficient decrease can be misleading in the context of Algorithm \ref{alg:linbregspecial} as it is not unusual for specific choices of $R$ that the function value of $E$ does not change for several iterations.

Our next result is a bound for the subgradients of the surrogate energy at the iterates computed with Algorithm \ref{alg:linbregspecial}. Note that the subdifferential of the surrogate objective reads as 
\begin{align*}
\partial F(x, y) = \left\{ \left( \left. \begin{array}{c}
\nabla E(x) + z_1 - y \\
z_2 - x
\end{array} \right) \, \right| \, z_1 \in \partial R(x), z_2 \in \partial R^\ast(y) \right\} \, ,
\end{align*}
which can for example be deduced from \cite{rockafellar2009variational}. With $q^{k + 1} \in \partial R(u^{k + 1})$, and the fact that $q^k \in \partial R(u^k)$ is equivalent to $u^k \in \partial R^\ast(q^k)$ (Lemma \ref{lem:convconj} in the appendix), we know that 
\begin{align}%
r^{k + 1} := \left( \begin{array}{c}
\nabla E(u^{k + 1}) + q^{k + 1} - q^k \\%
u^k - u^{k + 1}
\end{array} \right) \in \partial F(u^{k + 1}, q^k) = \partial F(s^{k + 1})\, .\label{eq:surrogatesubgrad}
\end{align}%
Subsequently, we want to show that the norm of this sequence of subgradients $\{ r^k \}_{k \in \N}$ is bounded by the iterates gap of the primal variable.

\begin{lemma}[A subgradient lower bound for the iterates gap]\label{thm:subgradientbound}
Let the same assumptions hold true as in Lemma \ref{thm:sufficientdecrease} and $\tau^k \geq \taumin := \inf_k \tau^k > 0$. Then the iterates of Algorithm \eqref{alg:linbregspecial} satisfy
\begin{align}
\| r^k \| \leq \rho_2 \| u^k - u^{k-1} \| \, , \label{eq:gradbounditeratesgap}
\end{align}
for $r^k \in \partial F(s^k)$ as defined in \eqref{eq:surrogatesubgrad}, $s^k := (u^k, q^{k - 1})$, $\rho_2 := \left(1 + L + 1 / \taumin \right)$ and $k \in \N$.
\begin{proof}
From \eqref{eq:surrogatesubgrad} we know 
\begin{align*}
\| r^k \| \leq \| \nabla E(u^k) + q^k - q^{k - 1} \| + \| u^k - u^{k - 1} \| \, .
\end{align*}
Together with \eqref{eq:linbregman:special:dual} we therefore estimate
\begin{align*}
\| r^k \| {} \leq {} &\left\| \nabla E(u^k) + q^k - q^{k-1} \right\| + \| u^k - u^{k - 1} \| \\
{} = {} &\left\| \nabla E(u^k) - \nabla E(u^{k-1}) + \frac{1}{\tau^{k-1}} \left( u^{k - 1} - u^k \right) \right\| + \| u^k - u^{k - 1} \| \, ,\\
{} \leq {} &\left(1 + L + \frac{1}{\taumin}\right) \| u^k - u^{k-1} \| = \rho_2 \| u^k - u^{k - 1} \| \, ,
\end{align*}
where we have made use of the Lipschitz-continuity of the gradient of $E$. 
\end{proof}
\end{lemma}

\begin{remark}
We want to point out that the Lipschitz-continuity of $\nabla E$ is not necessary if $R \equiv 0$. In that case it is easy to see that we can obtain the estimate
\begin{align*}
\| \nabla E(u^k) \| \leq \frac{1}{\taumin} \| u^{k + 1} - u^k \|
\end{align*}
instead of \eqref{eq:gradbounditeratesgap} (see also \cite{benning2016gradient}), without the use of Lipschitz-continuity. For the sufficient decrease Theorem \ref{thm:sufficientdecrease} it is already enough to choose $\tau^k$ such that $G := \frac{1}{2}\| \cdot \|^2 - \tau^k E$ is convex for all arguments and all $k \in \N$. This observation has already been made and exploited in \cite{bauschke2016descent,benning2016gradient,bolte2017first}. We also want to emphasise that the requirement of Lipschitz continuity can potentially be relaxed if backtracking strategies are incorporated into Algorithm \ref{alg:linbregspecial}.
\end{remark}

To conclude our convergence analysis we prove global convergence of Algorithm \ref{alg:linbregspecial} with the help of the Kurdyka-\L ojasiewicz (KL) property as defined in the appendix in Definition \ref{def:kl}. In order to apply the KL property, we have to verify some properties of the set of limit points. Let $\{ s^k \}_{k \in \N} = \{ (u^k, q^{k - 1}) \}_{k \in \N}$ be a sequence generated by Algorithm \ref{alg:linbregspecial} from starting points $u^0$ and $q^0$ with $q^0 \in \partial R(u^0)$. The set of limit points is defined as 

\begin{align*}
\omega(s^0) &:= \left\{ \overline{s} = (\overline{u}, \overline{q}) \in \R^n \times \R^n \, \left| \, \text{there exists an increasing sequence} \vphantom{\lim_{j \rightarrow \infty} s^{k_j} = \overline{s}} \right. \right. \\ 
&\qquad \left. \text{of integers $\{k_j \}_{j \in \N}$ such that $\lim_{j \rightarrow \infty} u^{k_j} = \overline{u}$ and $\lim_{j \rightarrow \infty} q^{k_j} = \overline{q}$} \right\} \, .
\end{align*}
Before we continue, we want to emphasise that the current assumptions on $E$ and $R$ are not sufficient in order to guarantee convergence of the dual variable, which we want to demonstrate with a simple counter example.

\begin{remark}\label{rem:critpoint}
Let $E(u) = (u + 1)^2 / 2$, and $R(u) = \chi_{\geq 0}(u)$ with
\begin{align*}
\chi_{\geq 0}(u) := \begin{cases} 0 & u \geq 0 \\ \infty & u < 0 \end{cases} \, .
\end{align*}
It is obvious that $E \in \FunProblem[1]$ and that the only critical point of $E$ is $\hat{u} = -1$. However, Algorithm \ref{alg:linbregspecial} can never converge to that point but will converge to $\overline{u} = 0$ due to the choice of $R$. This can be seen for instance for the choices $u^0 > 0$, $q^0 = 0$ and $\tau^k = 1$. Then the subsequent iterates are $u^k = 0$ and $q^k = u^0 - k$, thus, $u^k \rightarrow 0$ and $q^k \rightarrow -\infty$.
\end{remark}
For convex, quadratic fidelity terms (such as $E$ in the example above) it is sufficient to satisfy a source condition of the form $\partial R(\hat{u}) \neq \emptyset$ (which in Remark \ref{rem:critpoint} is clearly violated) in order to guarantee boundedness of the subgradients, see for instance \cite{frick2010regularization}. For general, non-convex terms $E$ it is not straight forward to adapt the concept of source conditions, which is why we are going to assume local boundedness of the subgradients instead.
\begin{definition}[Locally bounded subgradients]\label{def:boundsubgrad}
We say that $R$ has \emph{locally bounded subgradients} if for every compact set $U \subset \R^n$ there exists a constant $C \in (0, \infty)$ such that for all $v \in U$ and all $q \in \partial R(v)$ we have $\|q\| \leq C$.
\end{definition}
Boundedness is not a very restrictive requirement as it is for instance satisfied for the large class of Lipschitz-continuous functions. 
\begin{proposition}\label{prop:boundsubgrad}
Let $R \in \FunConvex$ be a (globally) Lipschitz continuous function in the sense of Definition \ref{def:lipschitz} in the appendix. Then $R$ has locally bounded subgradients.
\begin{proof}
From the convexity of $R$ we observe
\begin{align*}
\langle q, h \rangle \leq | R(v+h) - R(v) | \leq L \| h \| \, ,
\end{align*}
for any $U \subset \R^n$ and any $h, v \in U$ with $v + h \in U$ and $q \in \partial R(v)$. Taking the supremum over $h$ with $\|h\| \leq 1$ shows $\| q \| \leq L$, which proves the assertion.
\end{proof}
\end{proposition}
\begin{remark}
Note that every continuously differentiable function is already locally Lipschitz-continuous, and therefore has locally bounded gradients according to Proposition \ref{prop:boundsubgrad}.
\end{remark}



Before we show global convergence of Algorithm \ref{alg:linbregspecial} to a critical point of $E$, we need to verify that the surrogate function converges to $E$ on $\omega(s^0)$, that $\omega(s^0)$ is a non-empty, compact and connected set and that its primal limiting points form a subset of the set of critical points of $E$. The following lemma guarantees that for a sequence converging to a limit point we also know that the surrogate objective converges to the objective evaluated at this limit point.
\begin{lemma}\label{lem:critpoint}
Suppose $E \in \FunProblem$, $R \in \FunConvex$, and let $\overline{s} \in \omega(s^0)$. Then we already know 
\begin{align}
\lim_{k \rightarrow \infty} F(s^k) = F(\overline{s}) = E(\overline{u}) \, .\label{eq:surrogatelimit1}
\end{align}
\begin{proof}
Since $\overline{s}$ is a limit point of $\{ s^k \}_{k \in \N}$ we know that there exists a subsequence $\{ s^{k_j} \}_{j \in \N}$ with $\lim_{j \rightarrow \infty} s^{k_j} = \overline{s}$. Hence, we immediately obtain
\begin{align*}
\lim_{j \rightarrow \infty} F(s^{k_j}) &= \lim_{j \rightarrow \infty} \left\{ E(u^{k_j}) + D_R^{q^{k_j - 1}}(u^{k_j}, u^{k_j - 1}) \right\} = E(\overline{u}) \, ,
\end{align*}
due to the continuity of $E$ and $\lim_{j \rightarrow \infty} D_R^{q^{k_j - 1}}(u^{k_j}, u^{k_j - 1}) = 0$ as a result of Lemma \ref{thm:sufficientdecrease}. 
Since $\{ F(s^k) \}_{k \in \N}$ is also monotonically decreasing and bounded from below according to Remark \ref{rem:boundedsequence}, we can further conclude \eqref{eq:surrogatelimit1} as a consequence of the monotone convergence theorem.
\end{proof}
\end{lemma}

In addition to Lemma \ref{lem:critpoint}, the following lemma states that $\omega(s^0)$ is a non-empty, compact and connected set, and that the objective $F$ is constant on that set.
\begin{lemma}[{\cite[Lemma 5]{bolte2014proximal}}]\label{lem:limitpointsetproperties}
Suppose $E \in \FunProblem$ and that $R \in \FunConvex$ has locally bounded subgradients. Then the set $\omega(s^0)$ is a non-empty, compact and connected set, the surrogate objective $F$ is constant on $\omega(s^0)$ and we have $\lim_{k \rightarrow \infty} \dist(s^k, \omega(s^0)) = 0$.
\end{lemma}

We can further verify that the set of primal limiting points is a subset of the set of critical points of the energy $E$.
\begin{lemma}\label{lem:critpoint4}
Suppose $E \in \FunProblem$, and that $R \in \FunConvex$ has locally bounded subgradients. Then we have $\overline{u} \in \crit(E)$ for every $\overline{s} = (\overline{u}, \overline{q}) \in \omega(s^0)$. 
\begin{proof}
We prove this assertion by contradiction to the boundedness of the subgradients. Let $\overline{s} := (\overline{u}, \overline{q}) \in \omega(s^0)$, which means $\lim_{k\rightarrow \infty} u^k = \overline{u}$. Assume that $\nabla E(\overline{u}) \neq 0$ and let $c := \|\nabla E(\overline{u})\| > 0$. It follows from the subgradient update \eqref{eq:sumdualupdate} and the reverse triangle inequality $\left\|a + \sum_i a_i \right\| \geq \|a\| - \sum_i \|a_i\|$ that 
\begin{align*}
\|q^{k}\|
\geq \left\|\sum_{n = 0}^{k-1} \nabla E(\overline{u})\right\|- \|q^0\| - \sum_{n=0}^{k-1} \left[\frac{1}{\tau^n} \| u^{n+1} - u^n \| + \|\nabla E(u^n) - \nabla E(\overline{u})\|\right] \, .
\end{align*}
As $u^k \rightarrow \overline{u}$, there exists $K \in \N$ such that for all $n \geq K$ the bounds $\| u^n - \overline{u} \| \leq c \taumin / 8$ and $\|\nabla E(u^n) - \nabla E(\overline{u})\| \leq c/4$ hold. Thus, we have for all $n \geq K$ that $$1/\tau^n \|u^{n+1} - u^n\| + \|\nabla E(u^n) - \nabla E(\overline{u})\|\leq c/2 ,$$ and therefore
\begin{align*}
&\sum_{n=0}^{k-1} \left[\frac{1}{\tau^n} \| u^{n+1} - u^n \| + \|\nabla E(u^n) - \nabla E(\overline{u})\|\right] \\
&\leq \sum_{n=K}^{k-1} \left[\frac{1}{\tau^n} \| u^{n+1} - u^n \| + \|\nabla E(u^n) - \nabla E(\overline{u})\|\right] + \text{const}
\leq k c/2 + \text{const} ,
\end{align*}
for all $k \in \N$, with a constant independent of $k$. Combining these two estimates yields
\begin{align*}
\|q^k\|
&\geq \left\|\sum_{n = 0}^{k-1} \nabla E(\overline{u})\right\| - k c/2 + \text{const}
= k c/2 + \text{const} \, .
\end{align*}
Hence, we observe $\lim_{k\to\infty} \|q^k\| = \infty$, which is a contradiction to the boundedness of $\{q^k\}$. Thus, $\nabla E(\overline{u}) = 0$, which means $\overline{u} \in \crit(E)$.
\end{proof}
\end{lemma}

Now we have all the necessary ingredients to show the following global convergence result for Algorithm \ref{alg:linbregspecial}.

\begin{theorem}[Finite length property]\label{thm:finitelength1}
Suppose that $F$ is a KL function in the sense of Definition \ref{def:kl}. Further, assume $R \in \FunConvex$ with locally bounded subgradients. Let $\{ s^k \}_{k \in \N} = \{ (u^k, q^{k - 1} ) \}_{k \in \N}$ be a sequence generated by Algorithm \ref{alg:linbregspecial}. Then the sequence $\{ u^k \}_{k \in \N}$ has finite length, i.e.
\begin{align}
\sum_{k = 0}^\infty \| u^{k + 1} - u^k \| < \infty \, \text{.}\label{eq:finitelength}
\end{align}
\begin{proof}
We follow the steps of the proof of \cite[Theorem 1]{bolte2014proximal} but with non-trivial modifications. 

The sequence  $\{ u^k \}_{k \in \N}$ is bounded, which follows from the assumption $E \in \FunProblem$ and the monotonic decrease. Thus, we know that there exists a convergent subsequence $\{ u^{k_j} \}_{j \in \N}$ and $\overline{u} \in \R^n$ with $$\lim_{j \rightarrow \infty} u^{k_j} = \overline{u} \, .$$ 

As a consequence of Lemma \ref{lem:critpoint} we further know that $\lim_{k \rightarrow \infty} F(s^k) = F(\overline{s}) = E(\overline u)$. If there exists an index $l \in \N$ with $F(s^l) = E(\overline{u})$ the results follow trivially. If there does not exist such an index, we observe that for any $\eta > 0$ there exists an index $k_1$ such that 
$$E(\overline{u}) < F(s^k) < E(\overline{u}) + \eta$$
for all $k > k_1$. In addition, for any $\varepsilon > 0$ there exists an index $k_2$ with $$\dist(s^k, \omega(s^0)) < \varepsilon$$
for all $k > k_2$, due to Lemma \ref{lem:limitpointsetproperties}. Hence, if we choose $l := \max(k_1, k_2)$, we know that $u^k$ is in the set \eqref{eq:klset} for all $k > l$ according to Lemma \ref{lem:unifiedkl} in the appendix.

By Lemma \ref{lem:limitpointsetproperties}, $\omega(u^0)$ satisfies all the assumptions of Lemma \ref{lem:unifiedkl} and we have
\begin{align}
1 \leq \varphi^\prime( F(s^k) - E(\overline{u}) ) \, \dist(0, \partial F(s^k)) \label{eq:finitelength:main}
\end{align}
for all $k > l$. This inequality makes sense due to $F(s^k) > E(\overline{u})$ for all $k$.

\noindent From the concavity of $\varphi$ we know that
\begin{align*}
\varphi^\prime(x) \leq \frac{\varphi(x) - \varphi(y)}{x-y}
\end{align*}
holds for all $x, y\in [0,\eta), x > y$,
which we will use for the specific choices of $x = F(w^k) - E(\overline{u})$ and $y = F(s^{k + 1}) - E(\overline{u})$. Combining the latter with Lemma \ref{thm:sufficientdecrease} and abbreviating $$\varphi^k := \varphi( F(s^k) - E(\overline{u}) ) \,$$ yields
\begin{align}
\varphi^\prime(F(s^k) - E(\overline{u})) \leq \frac{\varphi^k - \varphi^{k+1}}{F(s^k) - F(s^{k+1})} \leq \frac{\varphi^k - \varphi^{k+1}}{\rho_1 \| u^{k + 1} - u^k \|^2} \, \text{.} \label{eq:kl:ineq}
\end{align}
Inserting \eqref{eq:kl:ineq} and the subgradient bound \eqref{eq:gradbounditeratesgap} into the KL inequality \eqref{eq:finitelength:main} leads to
\begin{align*}
\| u^{k+1} - u^k \|^2 \leq \frac{\rho_2}{\rho_1} (\varphi^k - \varphi^{k+1}) \| u^k - u^{k-1} \| \, .
\end{align*}

Taking the square root, multiplying by 2 and using Young's inequality of the form $2 \sqrt{ab} \leq a + b$ then yields
\begin{align*}
2 \| u^{k+1} - u^k \| \leq \frac{\rho_2}{\rho_1} (\varphi^k - \varphi^{k+1}) + \| u^k - u^{k-1} \| \, .
\end{align*}
Subtracting $\| u^{k+1} - u^k \|$ and summing from $k=l, \ldots, N$ leads to
\begin{align*}
\sum_{k = l}^{N} \| u^{k + 1} - u^k \| 
&\leq \frac{\rho_2}{\rho_1} (\varphi^l - \varphi^{N+1}) + \| u^{l} - u^{l-1} \| - \| u^{N+1} - u^N \| \\
&\leq \frac{\rho_2}{\rho_1} \varphi^l + \| u^{l} - u^{l-1} \| < \infty \, ,
\end{align*}
and hence, we obtain the finite length property by taking the limit $N \rightarrow \infty$.
\end{proof}
\end{theorem}

\begin{corollary}[Convergence]\label{cor:convergence}
Under the same assumptions as Theorem \ref{thm:finitelength1}, the sequence $\{ u^k \}_{k \in \N}$ converges to a critical point of $E$.
\begin{proof}
As in the proof of \cite[Theorem 1 (ii)]{bolte2014proximal}, the finite length property Theorem \ref{thm:finitelength1} implies $\sum_{k = l}^{N} \| u^{k+1} - u^k\| \rightarrow 0$ for $N \rightarrow \infty$. Thus, for any $s \geq r \geq l$ we have 
\begin{align*}
\| u^s - u^r \| = \left\| \sum_{k = r}^{s - 1} u^{k + 1} - u^k \right\| \leq \sum_{k = r}^{s - 1} \| u^{k + 1} - u^k \| \leq \sum_{k = l}^{\infty} \| u^{k+1} - u^k\| \, .
\end{align*}
This shows that $\{ u^k \}_{k \in \N}$ is a Cauchy sequence and, thus, is convergent. According to Lemma \ref{lem:critpoint4} its limit is a critical point of $E$.
\end{proof}
\end{corollary}

\subsection{Global convergence in the absence of locally bounded subgradients}\label{subsec:noboundsubgrad}
In the previous section we have made the assumption that the subgradients of $R$ have to be locally bounded in order to guarantee convergence of the primal iterates to a critical point of $E$. In Remark \ref{rem:critpoint} we have seen an example for which the subgradients of $R$ diverge, but the primal iterates still converge, just not to a critical point of $E$. This leaves us with two open questions: 1) could we prove convergence of the primal iterates without boundedness of the dual iterates and 2) would the limit (if it exists) be a critical point of some other energy? 
It might be possible to answer the first question by slightly modifying Definition \ref{def:kl} and Lemma \ref{lem:unifiedkl} in the appendix, as well as Lemma \ref{lem:limitpointsetproperties} to accommodate the fact that the surrogate function is also constant on the set of limiting points that only depends on the primal variable (which we denote by $\omega(u^0)$ for convenience). A potential modification of \eqref{eq:klset} in Lemma \ref{lem:limitpointsetproperties} could for instance be
\begin{align*}
    \{ u, q \in \R^n \, | \, \dist(u, \omega(u^0)) < \varepsilon \} \cap \{ u, q \in \R^n \, | \, E(\overline{u}) \leq F(u, q) \leq E(\overline{u}) + \eta \} \, ,
\end{align*}
where $\overline{u} \in \omega(u^0)$. Note that this modification would not affect the finite length proof of Theorem \ref{thm:finitelength1} and therefore would still imply global convergence, but not necessarily to a critical point of $E$. 
Remark \ref{rem:critpoint} leaves room for speculation whether an answer to the second question is that the primal iterates converge to a critical point of $E + \chi_{\dom(R)}$, where $\chi_{\dom(R)}$ denotes the characteristic function over the effective domain of $R$. Proving this, however, is beyond the scope of this paper.

\changed{\subsection{Limitations of the convergence analysis and possible remedies}\label{subsec:limitations}
The convergence analysis presented in this paper relies on the fact that the function $E$ satisfies $E \in \FunSmooth{L}$, which is often restrictive for practical applications. Even simple functions such as the blind deconvolution data fidelity term from Section \ref{sec:motivation} are not globally $L$-smooth. Remedies are the use of an alternating version of Algorithm \ref{alg:linbregspecial} in the spirit of \cite{bolte2014proximal} and to make use of local smoothness of the functions with fixed variables. \added{For two variables $u_1$ and $u_2$, such a scheme is of the form
\begin{align*}
u_1^{k + 1} &= \argmin_{u_1 \in \uc} \left\{ \tau^k_1 \langle \nabla_1 E(u_1^k, u_2^k), u_1 - u_1^k\rangle + D_{J_1}^{p_1^k}(u_1, u_1^k) \right\} \\
p_1^{k +1} &= p_1^k - \tau^k_1 \nabla_1 E(u_1^k, u_2^k) \, ,\\
u_2^{k + 1} &= \argmin_{u_2 \in \uc} \left\{ \tau^k_2 \langle \nabla_2 E(u_1^{k + 1}, u_2^k), u_2 - u_2^k\rangle + D_{J_2}^{p_2^k}(u_2, u_2^k) \right\} \\
p_2^{k +1} &= p_2^k - \tau^k_2 \nabla_2 E(u_1^{k + 1}, u_2^k) \, ,
\end{align*}
assuming a separable structure of $J(u_1, u_2) = J_1(u_1) + J_2(u_2)$. Here $\nabla_1$ and $\nabla_2$ refer to the partial gradients of $E$ with respect to $u_1$ and $u_2$, and $p_1^k \in \partial J_1(u_1^k)$ and $p_2^k \in \partial J_2(u_2^k)$ are subgradients in the subdifferential of $J_1$ and $J_2$, respectively.} The analysis of such a scheme should be relatively straight-forward, but is beyond the scope of this work. 

Another limitation in terms of convergence analysis that becomes obvious from the motivating example in Section \ref{sec:motivation} is the use of characteristic functions. If we incorporate them in the function $R$, we run into the issues outlined in Section \ref{subsec:noboundsubgrad}. If we add them to the objective function $E$, we lose the continuity and differentiability. A remedy for the blind deconvolution example (and many similar examples) in Section \ref{sec:motivation} is that for the convolution kernel the additional Bregman function $R$ is simply zero, so that the algorithm merely has to perform a proximal point step in the direction of the convolution kernel. The convergence analysis in such a setting is straight-forward, but we did not include it in order not to complicate notation. Alternatively, one could replace the characteristic function with its Moreau--Yosida envelope.}

This concludes the theoretical analysis of Algorithm \ref{alg:linbregspecial}. In the following two sections we are going to discuss three applications, their mathematical modelling in the context of Algorithm \ref{alg:linbregspecial} and their numerical results.

\section{Applications}\label{sec:apps}
We demonstrate the capabilities of the linearised Bregman iteration by using it to approximately minimise several non-convex minimisation problems. We say approximately, as we do not exactly minimise the corresponding objective functions, but rather compute iteratively regularised solutions to the associated inverse problems via early stopping of the iteration.

\subsection{Parallel Magnetic Resonance Imaging}\label{subsec:pmri}
In (standard) Magnetic Resonance Imaging (MRI) the goal is to recover the spin-proton density from sub-sampled Fourier measurements that were obtained with a single radio-frequency (RF) coil. In parallel MRI, multiple RF coils are used for taking measurements, thus allowing to recover the spin-proton density from more measurements compared to the standard case. This, however, comes at the cost of having to model the sensitivities of the individual RF coils w.r.t. the measured material. We basically follow the mathematical modelling of \cite{Pruessmann1999, uecker2008image} and describe the recovery of the spin-proton density and the RF coil sensitivities as the minimisation of the following energy function:
\begin{align}
E(u, b_1,\! \ldots\!, b_s) \! := \! \frac{1}{2} \! \sum_{j = 1}^s \| \mathcal{S}(\mathcal{F}((K(u, b_1,\! \ldots\!, b_s))_j )) \! - \! f_j \|_2^2 \! + \! \frac{\epsilon}{2} \! \left( \! \| u \|^2 + \sum_{j=1}^s \| b_j \|^2 \! \right) \! .\label{eq:parallelmri}
\end{align}
Here $\mathcal{F} \in \C^{n \times n}$ is the (discrete) Fourier transform, $\mathcal{S} \in \{0, 1\}^{m \times n}$ is a sub-sampling operator, $K$ is the non-linear operator $K(u, b_1, \ldots, b_s) = (u b_1, u b_2, \ldots, u b_s )^T$, $u$ denotes the spin-proton density, $b_1, b_2, \ldots, b_s$ the $s$ coil sensitivities, $f_1, \ldots, f_s$ the corresponding sub-sampled k-space data and $\epsilon > 0$ is a scalar parameter that ensures bounded level-sets of $E$. Since $\C$ has the same topology as $\R \times \R$, we can formally treat all variables as variables in $\R^{2n}$. \changed{Note that $E$ as defined in \eqref{eq:parallelmri} is not globally $L$-smooth, which is why we also assume that we choose parameters and initial values such that our sequence $\{ u^k \}_{k \in \N}$ of primal variables generated by Algorithm \ref{alg:linbregspecial} satisfies
\begin{align*}
\| \nabla E(u^k) - \nabla E(u^{k - 1}) \|_2 \leq L^k \| u^k - u^{k - 1} \|_2 \, , 
\end{align*}
for a sequence $\{ L^k \}_{k \in \N}$ of positive constants. Hence, $E \in \FunProblem[L^k]$, which means that $E$ is (locally) $L^k$-smooth, respectively $\nabla E$ is (locally) $L^k$-Lipschitz-continuous in the sense of Definition \ref{def:lipschitz}. Furthermore, we assume that the sequence $\{ L^k \}_{k \in \N}$ is bounded from above, i.e. $L^k \leq L$ for all $k \in \N$, and consequently $E \in \FunProblem[L]$. It is not necessarily straight-forward to prove existence of $L$ a-priori, but it is relatively easy to validate it a-posteriori. Note that, alternatively, one could use an alternating version of Algorithm \ref{alg:linbregspecial} as discussed in Section \ref{subsec:limitations}.}

The inverse problem of parallel MRI has been subject in numerous research publications \cite{ramani2011parallel,knoll2012parallel,benning2015preconditioned}. We follow a different methodology here and apply Algorithm \ref{alg:linbregspecial} to approximately minimise \eqref{eq:parallelmri} with the following configuration. We choose the function $R$ to be of the form
\begin{align*}
R(u, b_1, \ldots, b_s) &= R_1(u) + \sum_{j = 1}^s R_2(b_j) \, ,
\intertext{with}
R_1(u) &= \alpha_0 \text{TV}(u) = \alpha_0 \| | \nabla u | \|_1
\intertext{and}
R_2(b_j) &= \alpha_j \sum_{l = 1}^n w_l \left| (C \, b_j)_l \right|, \qquad \forall j \in \{ 1, \ldots, s \} \, .
\end{align*}
Here $\nabla$ denotes a discrete finite forward difference approximation of the gradient, $| \cdot |$ is the Euclidean vector norm, $C$ denotes the discrete two-dimensional cosine transform, $\{ w_l \}_{l \in \{1, \ldots, n\}}$ is a set of weighting-coefficients and $\alpha_0, \ldots, \alpha_{s + 1}$ are positive scaling parameters. Note that all functions are chosen to be semi-algebraic, and semi-algebraic functions and their additive compositions are KL functions (see \cite{attouch2009convergence,attouch2010proximal,attouch2013convergence}). Iterating Algorithm \ref{alg:linbregspecial} for too long may lead to unstable minimisers of \eqref{eq:parallelmri} in case the k-space data $f_1, \ldots, f_s$ are noisy, which is why we are going to apply Morozov's discrepancy principle \cite{morozov2012methods} as a stopping criterion to stop the iteration early (see also \cite{osher2005iterative,garrigos2016iterative,matet2017don}, and \cite{schopfer2006nonlinear,Bachmayr2009,kaltenbacher2009iterative} in the context of nonlinear inverse problems), i.e. we stop the iteration as soon as
\begin{align}
E(u, b_1, \ldots, b_s) \leq \eta \label{eq:morozov}
\end{align}
is satisfied, for some $\eta > 0$. Usually $\eta$ depends on the variance of the normal-distributed noise. 

\subsection{Blind deconvolution}
Blind deconvolution is extensively discussed in the literature, e.g. \cite{Kundur1996blinddeblurring,Chan2005book,Campisi2016book} and the references therein, with several approaches for which the convergence proofs also rely on the KL inequality \cite{bolte2010alternating,repetti2014euclid,chouzenoux2016block}. We follow the same setting as in Section \ref{sec:motivation} (with additional regularisation as in \eqref{eq:parallelmri} in order to guarantee bounded level-sets) and make the assumptions that the blur-free image $u$ has low total variation and that the kernel $h$ satisfies a simplex constraint, i.e. all entries are non-negative and sum up to one. The assumption of low total variation can for instance be motivated by \cite{chan1998total}, but as as we have seen in Section \ref{sec:motivation}, minimising $E$ with some additional total variation regularisation does often not lead to visually satisfactory results. We therefore apply Algorithm \ref{alg:linbregspecial} with $R :\R^n \times \R^r \rightarrow \R$ defined as
\begin{align*}
R(u, h) = \alpha \text{TV}(u) \, ,
\end{align*}
for $\alpha \geq 0$. All functions are semi-algebraic, and \changed{we make the same local smoothness assumption as in Section \ref{subsec:pmri}}. In case of noisy data, we will proceed as in Section \ref{subsec:pmri} and stop the iteration via the discrepancy principle.

\subsection{Classification}\label{subseq:classification}
The last application that we want to discuss is the classification of images. Given a set $D \in \R^{s \times r}$ of $r$ training images (with $s$ pixel each) in column vector form, we want to train a neural network to classify those images. We do so by learning the parameters $(A_1, \ldots, A_l)$ of the $l$-layer neural network
\begin{align*}
\rho(x) := \rho_1( A_1 \rho_2( A_2 \ldots \rho_l( A_l x ) ) \ldots )
\end{align*}
in a supervised fashion. Here the parameters $A_j \in \R^{m_j \times n_j}$ are matrices of different size, and the functions $\{ \rho_j \}_{j = 1}^l$ are so-called activity functions of the neural net. Typical choices for activity functions are $\max$- and $\min$-functions, also known as rectifier. However, due to their non-differentiability it is common to approximate them with either the pointwise smooth-$\max$-function, i.e.
\begin{align*}
\rho_j(x, c, \beta) := \frac{x \exp(\beta x) + c \exp(\beta c)}{\exp(\beta x) + \exp(\beta c)} \, ,
\end{align*} 
for $x \in \R$ and constants $\beta,  \in \R$, or the soft-$\max$-function, i.e.
\begin{align*}
\rho_j(x)_i = \frac{\exp( x_i )}{\sum_{l = 1}^m \exp( x_l ) } \, ,
\end{align*}
for $x \in \R^m$. The latter has the advantage that the function output automatically satisfies the simplex constraint. 

\begin{figure}[!t]
\begin{center}
\subfloat[Fully sampled]{\includegraphics[width=0.32\textwidth,trim={2.8cm 0.75cm 2cm 0.75cm},clip]{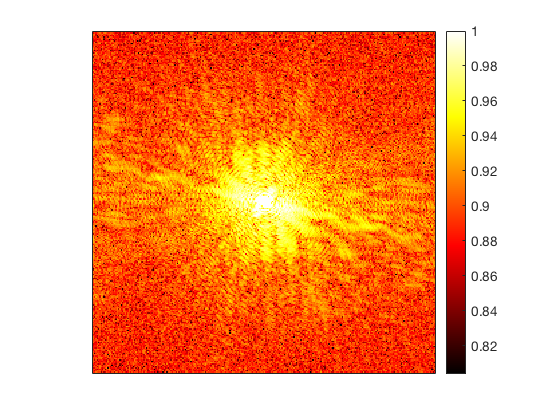}\label{subfig:fullysampledkspace}}
\subfloat[Recon. from \ref{subfig:fullysampledkspace}]{\includegraphics[width=0.32\textwidth,trim={2.8cm 0.75cm 2cm 0.75cm},clip]{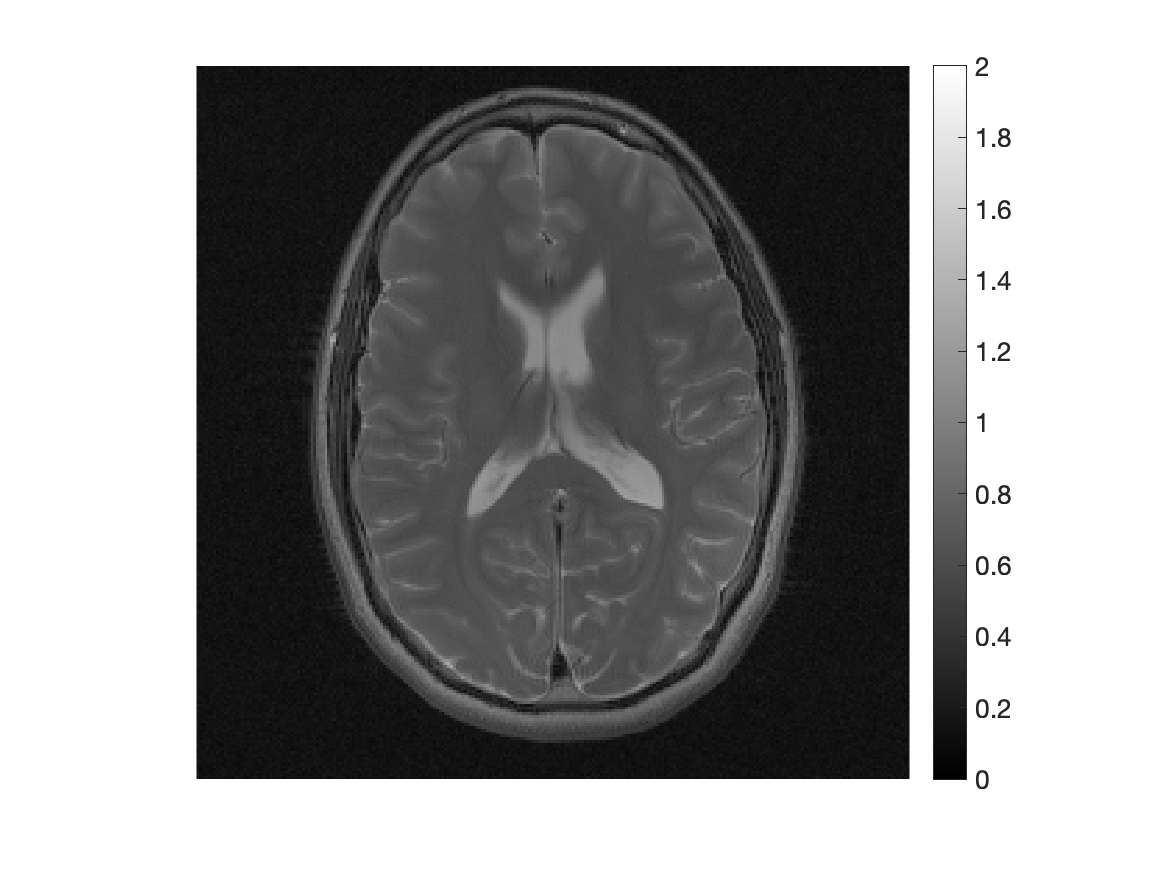}\label{subfig:fullysampledgrad-descent}}
\subfloat[Recon. from \ref{subfig:fullysampledkspace}]{\includegraphics[width=0.32\textwidth,trim={2.8cm 0.75cm 2cm 0.75cm},clip]{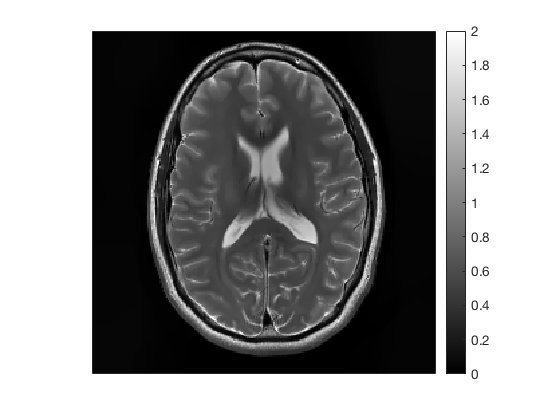}\label{subfig:fullysampledreconspin}}\\
\subfloat[Subsampled]{\includegraphics[width=0.32\textwidth,trim={2.8cm 0.75cm 2cm 0.75cm},clip]{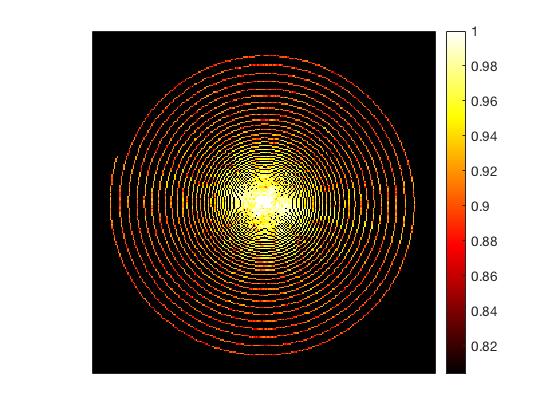}\label{subfig:subsampkspace}}
\subfloat[Recon. from \ref{subfig:subsampkspace}]{\includegraphics[width=0.32\textwidth,trim={2.8cm 0.75cm 2cm 0.75cm},clip]{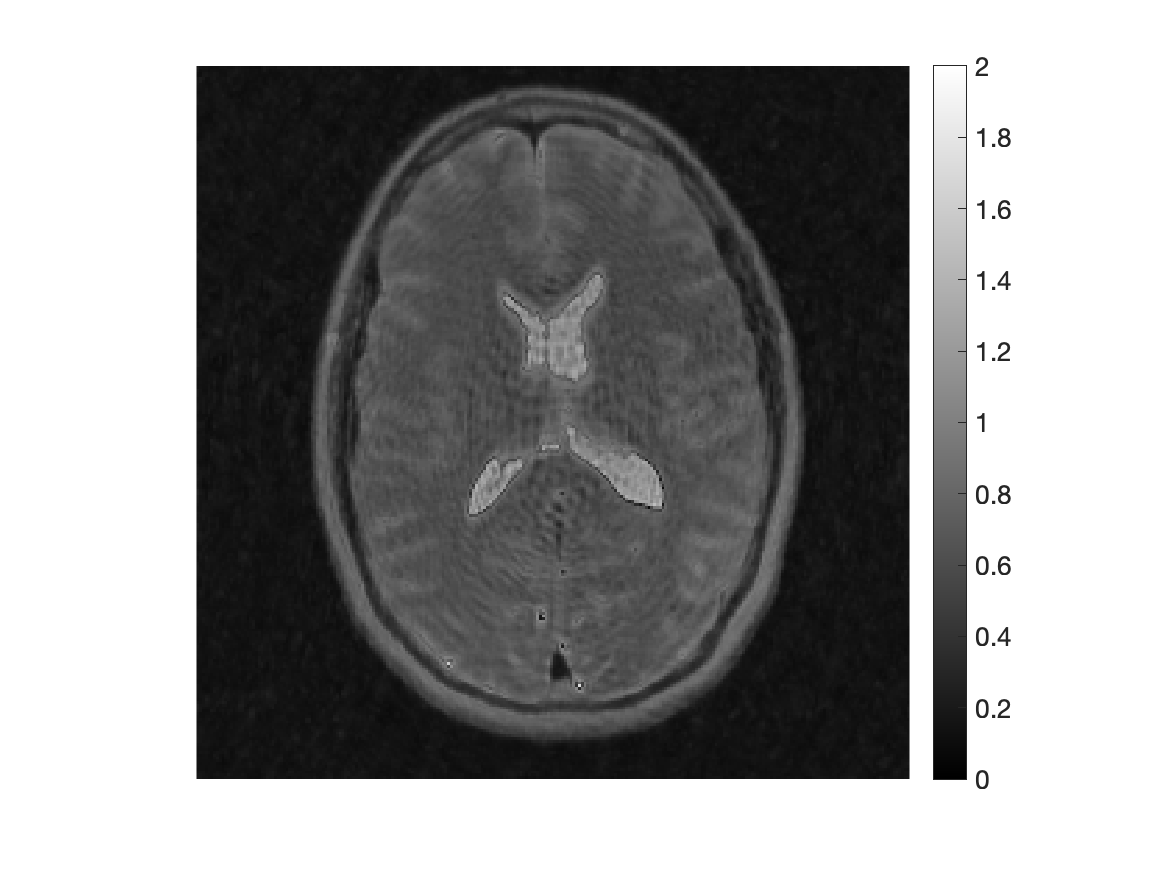}\label{subfig:subsampledgrad-descent}}
\subfloat[Recon. from \ref{subfig:subsampkspace}]{\includegraphics[width=0.32\textwidth,trim={2.8cm 0.75cm 2cm 0.75cm},clip]{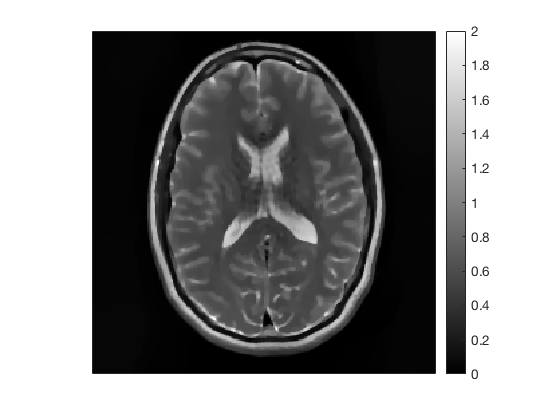}\label{subfig:subsampledreconspin}}\\
\subfloat[Convergence plot from \ref{subfig:fullysampledkspace}]{ \includegraphics[width=0.35\textwidth]{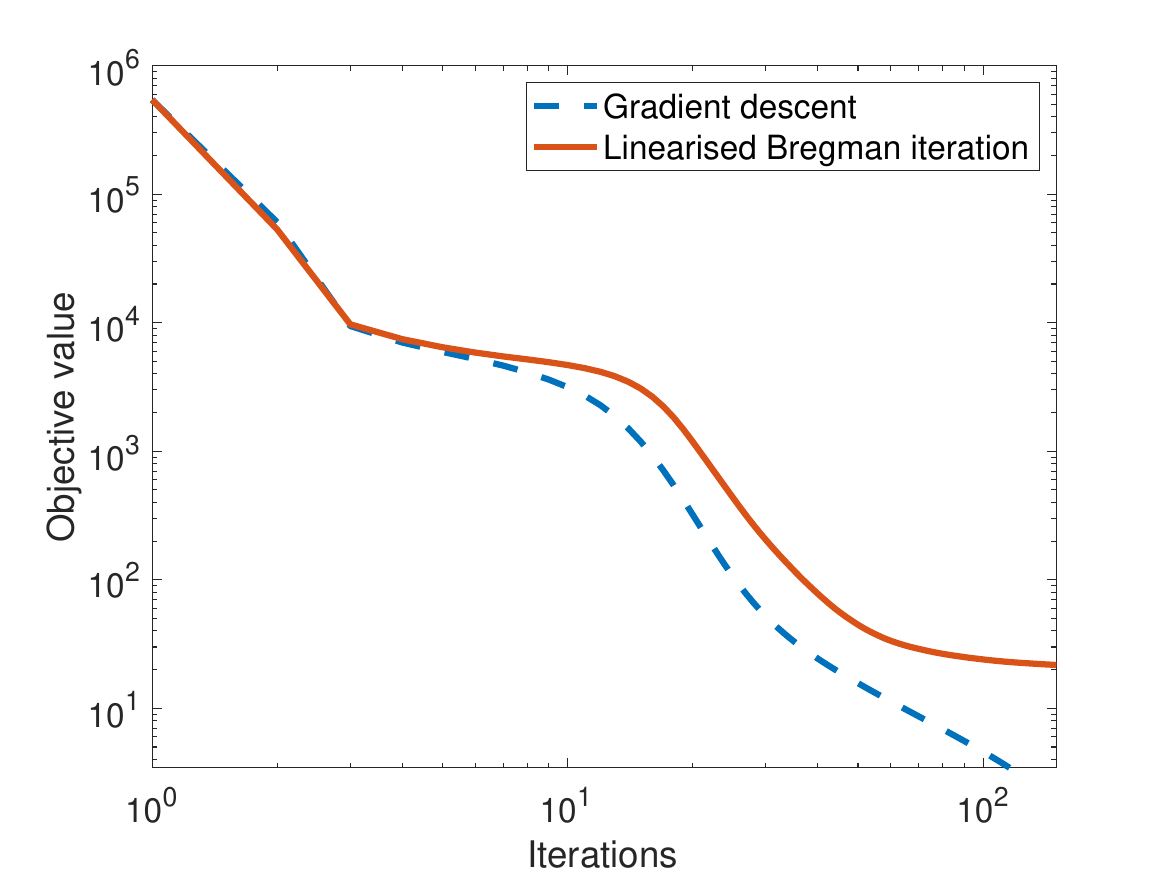}\label{subfig:fullysampled-convergence}}
\subfloat[Convergence plot from \ref{subfig:subsampkspace}]{\includegraphics[width=0.35\textwidth]{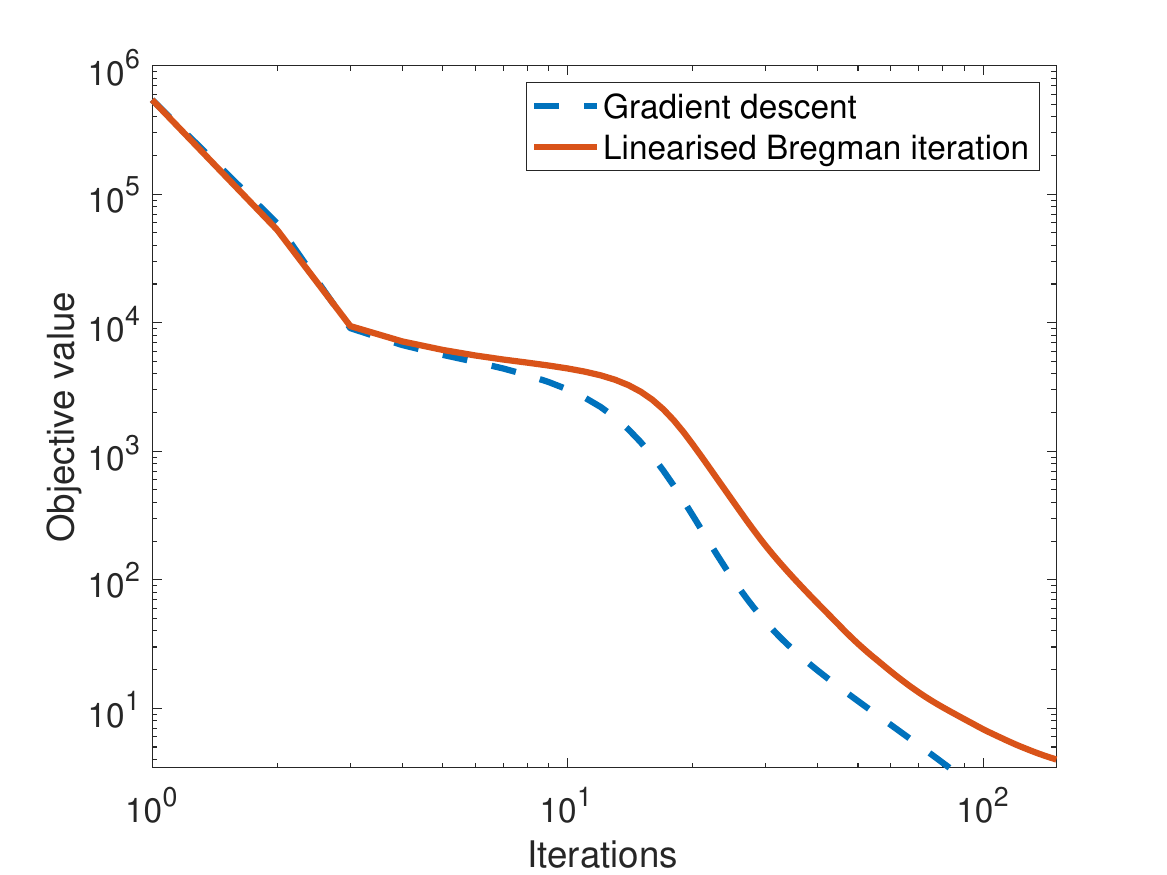}\label{subfig:subsampled-convergence}}
\end{center}
\caption{Figure \ref{subfig:fullysampledkspace} shows a log-plot of the modulus of the fully sampled k-space data of the first coil taken from \cite{knoll2012parallel}. Figure \ref{subfig:fullysampledgrad-descent} shows the reconstruction of the spin proton density from the data visualised in Figure \ref{subfig:fullysampledkspace} via gradient descent, whereas Figure \ref{subfig:fullysampledreconspin} shows the reconstruction of the spin proton density from the same data but via Algorithm \ref{alg:linbregspecial}. In Figure \ref{subfig:subsampkspace} we see roughly 25 \% of the k-space data visualised in Figure \ref{subfig:fullysampledkspace}, sampled on a spiral on a cartesian grid \cite{benning2014phase}. Figure \ref{subfig:subsampledgrad-descent} shows the reconstruction of the spin proton density from this subsampled k-space data with gradient descent, while Figure \ref{subfig:subsampledreconspin} shows the reconstruction of the spin proton density from the same data but with Algorithm \ref{alg:linbregspecial}. Figure \ref{subfig:fullysampled-convergence} and Figure \ref{subfig:subsampled-convergence} are showing the convergence plots in terms of energy decrease per iteration for the reconstructions that are obtained from the k-space data shown in Figure \ref{subfig:fullysampledkspace}, respectively in Figure \ref{subfig:subsampkspace}.}
\end{figure}

Note that if each function $\rho_j(A_j x)$ is chosen to be semi-algebraic, the composition $\rho$ is also semi-algebraic, see \cite[Proposition 2.2.10]{aizenbud2008schwartz}. If we choose $\rho_j(y) := \min(1, \max(0, y) )$ for all $j \in \{1, \ldots, l\}$ for instance, we can then show that also $\rho$ is semi-algebraic.

\noindent Defining a nonlinear operator $K(A_1, A_2, \ldots, A_l) := \rho_1( A_1 \rho_2( A_2 \ldots \rho_l( A_l D ) ) \ldots )$ for a given matrix $D$ and a given label matrix $Y \in \R^{m_1 \times r}$, we aim to minimise 
\begin{align}
E(A_1, A_2, \ldots, A_l) := \mathcal{D}(K(A_1, A_2, \ldots, A_l), Y) + \frac{\epsilon}{2} \sum_{j = 1}^l \| A_j \|_{\text{Fro}}^2 \, ,\label{eq:learning}
\end{align}
where $\mathcal{D}:\R^{m_1 \times r} \times \R^{m_1 \times r} \rightarrow \R$ denotes a function that measures the distance between its arguments in some sense. Our choice for $\mathcal{D}$ is simply the squared Frobenius norm $\mathcal{D}(X, Y) = \frac{1}{2}\| X - Y \|_{\text{Fro}}^2$ but other choices are possible. As mentioned earlier, the whole objective $E$ can be made a KL function, if for instance $\mathcal{D}$ and $\rho$ are chosen to be semi-algebraic, as their composition will also be semi-algebraic.

As in the previous sections, we aim to minimise \eqref{eq:learning} with Algorithm \ref{alg:linbregspecial} \changed{and make the same local smoothness assumption as before}. This time we choose $R(A_1, \ldots, A_l) = \sum_{j = 1}^l \alpha_j \| A_j \|_{\ast}$. Here $\{ \alpha_j \}_{j}^l$ is a set of positive scaling parameters, and $\| X \|_\ast := \sum_{i = 1}^{\text{rank}(X)} \sigma_i $ is the one norm of the singular values $\{ \sigma_i \}_{i = 1}^{\text{rank}(X)}$ of the argument $X$, also known as the nuclear norm. The rationale behind this choice for $R$ is that we can create iterates where the ranks of the individual matrices are steadily increasing. This way we control the number of effective parameters and do not fit all parameters right from the start.

\begin{figure}[!t]
\begin{center}
\subfloat[First coil]{\includegraphics[width=\coillength\textwidth]{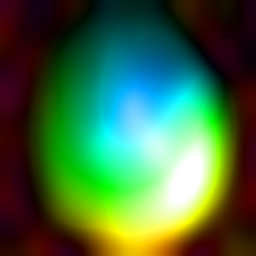}\label{subfig:fullysampledcoil1}}\vspace{0.05cm}
\subfloat[Second coil]{\includegraphics[width=\coillength\textwidth]{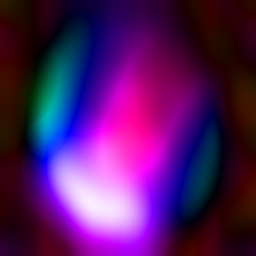}}\vspace{0.05cm}
\subfloat[Third coil]{\includegraphics[width=\coillength\textwidth]{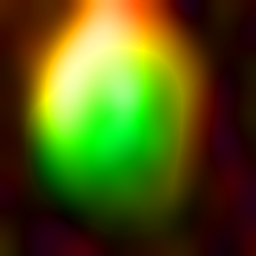}}\vspace{0.05cm}
\subfloat[Fourth coil]{\includegraphics[width=\coillength\textwidth]{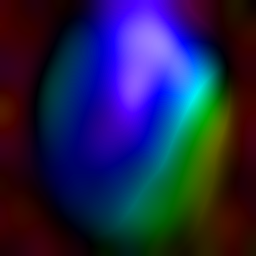}\label{subfig:fullysampledcoil4}}\vspace{0.05cm}
\phantom{\DrawColorwheel}\\
\subfloat[First coil]{\includegraphics[width=\coillength\textwidth]{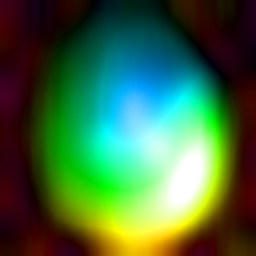}\label{subfig:subsampledcoil1}}\vspace{0.05cm}
\subfloat[Second coil]{\includegraphics[width=\coillength\textwidth]{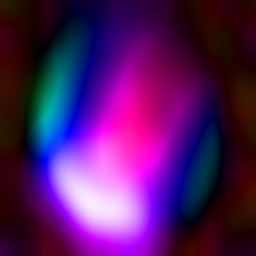}}\vspace{0.05cm}
\subfloat[Third coil]{\includegraphics[width=\coillength\textwidth]{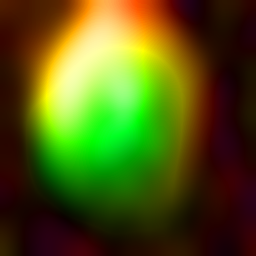}}\vspace{0.05cm}
\subfloat[Fourth coil]{\includegraphics[width=\coillength\textwidth]{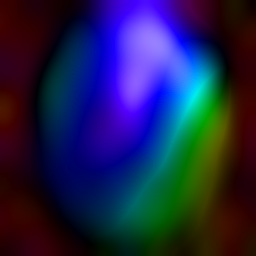}\label{subfig:subsampledcoil4}}\vspace{0.05cm}
\DrawColorwheel%
\end{center}
\caption{Figure \ref{subfig:fullysampledcoil1} - \ref{subfig:fullysampledcoil4} show the reconstructions of the coil sensitivities from the fully sampled data. Figure \ref{subfig:subsampledcoil1} - \ref{subfig:subsampledcoil4} show the reconstructions of the same quantities from the sub-sampled data.}
\end{figure}

\section{Numerical Results}\label{sec:numresults}
We demonstrate the particular properties and idiosyncrasies of Algorithm \ref{alg:linbregspecial} by computing several numerical solutions to the problems described in Section \ref{sec:apps}. All results have been computed with MATLAB R2017b. The code for the following examples is available at \url{https://doi.org/10.17863/CAM.16931} and can be used under the Creative Commons Attribution (CC BY) license once the article is accepted for publication.

Notably, all regularisation parameters that ensure boundedness of the level-sets are set to the smallest possible value ($\epsilon =$ machine accuracy) in practice. Since we do not use explicit Lipschitz constants, we employ a na\"{i}ve backtracking strategy for the variable stepsize $\{ \tau^k \}_{k \in \N}$. We start with an initial stepsize $\tau^0 > 0$ and check after each iteration whether $E(u^{k + 1}) \leq E(u^k) + \varepsilon$ is satisfied. Here, $\varepsilon > 0$ is a small constant that accounts for numerical rounding errors that may cause $E(u^{k + 1}) > E(u^k)$ when $E(u^{k + 1}) \approx E(u^k)$. If the decrease is satisfied, we set $\tau^{k + 1} = \tau^k$; otherwise we set $\tau^{k + 1} = (3 \tau^k)/4$ and backtrack again until we get a decrease. We want to emphasise that more sophisticated backtracking approaches can be used; we found, however, that the na\"{i}ve strategy that we use already works well for the computational results shown in the following subsections.

\subsection{Parallel MRI}
We compute parallel MRI reconstructions from real k-space data. We use data from a T2-weighted TSE scan of a transaxial slice of a brain acquired with a four-channel head-coil in \cite{knoll_florian_2010_800525}. A reconstruction from fully sampled data is taken as a ground truth. The spiral sub-sampling is simulated by point-wise multiplication of the k-space data with the spiral pattern visualised in Figure \ref{subfig:subsampkspace}. We initialise with $u^0 = 2 \times \textbf{1}^{65536 \times 1}$ and $b_j^0 = \textbf{1}^{65536 \times 1}$ for $j \in \{1, \ldots, 4\}$, and compute a $q^0 \in \partial R(u^0)$.

With the parameters $\alpha_j = 1$ for all $j \in \{0, \ldots, s\}$, $\tau^0 = 1/2$, $w_1 = w_2 = w_{\sqrt{n} + 1} = w_{\sqrt{n} + 2} = 10^{-6}$ and $w_l = 5$ for $l \in \{1, \ldots, n\}\setminus \{ 1, 2, \sqrt{n} + 1, \sqrt{n} + 2\}$, and $\eta = 3.45$ we obtain the spin proton density reconstruction visualised in Figure \ref{subfig:fullysampledreconspin}, as well as the coil sensitivity reconstructions in Figure \ref{subfig:fullysampledcoil1} - \ref{subfig:fullysampledcoil4}. In Figure \ref{subfig:subsampledreconspin} and Figure \ref{subfig:subsampledcoil1} - \ref{subfig:subsampledcoil4} we show the results of the reconstructions from sub-sampled data using the sub-sampling scheme in Figure \ref{subfig:subsampkspace}.

\begin{figure}[!t]\label{fig:blinddeconvolution}
\begin{center}
\subfloat[Original image]{\includegraphics[width=0.49\textwidth]{Images/Pixel.png}\label{subfig:deconvolution1}}\vspace{0.05cm}
\subfloat[Noisy, blurred image]{\includegraphics[width=0.49\textwidth]{Images/Pixelblurred.png}\label{subfig:deconvolution2}}\\
\subfloat[$\alpha = 0.1$]{\includegraphics[width=0.49\textwidth]{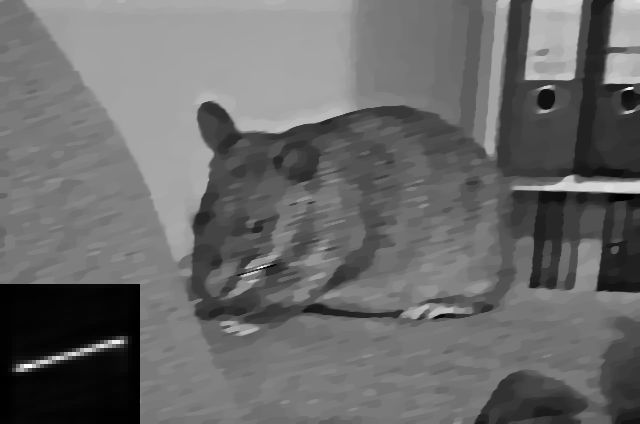}\label{subfig:blinddeconvolution1}}\vspace{0.05cm}
\subfloat[$\alpha = 10^{-3}$]{\includegraphics[width=0.49\textwidth]{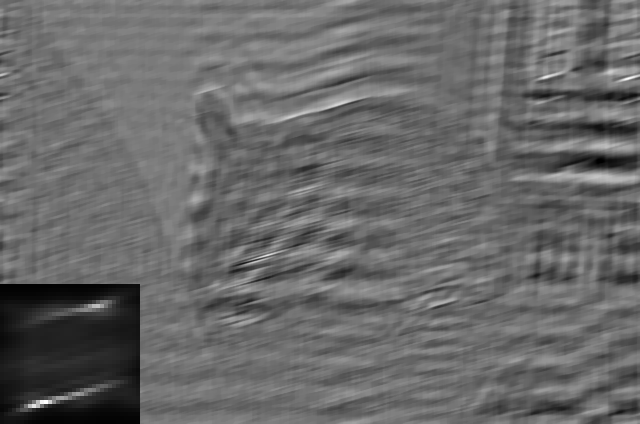}\label{subfig:blinddeconvolution2}}
\end{center}
\caption{Figure \ref{subfig:deconvolution1} shows an image of Pixel the Gambian pouched rat. Figure \ref{subfig:deconvolution2} shows a motion-blurred version of that image, together with some added normal distributed noise. The corresponding convolution kernel is depicted in the bottom left corner. Figure \ref{subfig:blinddeconvolution1} visualises the reconstruction of the image and the convolution kernel with Algorithm \ref{alg:linbregspecial} for the choice $\alpha = 10^{-1}$. Figure \ref{subfig:blinddeconvolution2} show the reconstructions of the same quantities for the choice $\alpha = 10^{-3}$. We clearly see that a larger choice of $\alpha$ results in a regular solution, whereas a smaller $\alpha$ will mimic traditional gradient descent with almost no additional regularity of the reconstruction.} 
\end{figure}

\subsection{Blind deconvolution}
To simulate blurring of a gray-scale image $f_{\rm orig} \in \R^{424 \times 640}$ we
subtract its mean, normalise it and subsequently blur $f_{\rm orig}$ with a motion-blur filter $h \in \R^{9 \times 31}$. The filter was obtained with the MATLAB\textcopyright-command \verb|fspecial('motion', 30, 15)|, and we assume periodic boundary conditions for the blurring process. Subsequently we add normally distributed noise with mean zero and standard deviation $\sigma = 10^{-4}$ to obtain a blurry and noisy image $f$ with ground truth $f_{\rm orig}$. Both $f_{\rm orig}$ and $f$, as well as $h$ are visualised in Figure \ref{fig:blinddeconvolution}. 

We use $f$ as our input image for Algorithm \ref{alg:linbregspecial}. We initialise Algorithm \ref{alg:linbregspecial} with $u^0 = 0$ and $q^0 = 0$. We choose $h_0 = 1/(r^2) \times \textbf{1}_{r \times r}$ for $r = 35$ to ensure that $h_0$ satisfies the simplex constraint. We set $\tau^0 = 2$ and pick $\alpha \in \{10^{-1}, 10^{-3}\}$. We then iterate Algorithm \ref{alg:linbregspecial} until the discrepancy principle is violated for $\eta = (1.2 \sigma^2)/(2 \sqrt{424 \times 640})$. The inner total variation sub-problem is solved with the primal-dual hybrid gradient method \cite{zhu2008efficient,pock2009algorithm,esser2010general,chambolle2011first,chambolle2016introduction}. The results are visualised in Figure \ref{fig:blinddeconvolution}.

\begin{figure}[!t]\label{fig:classification}
\begin{center}
\subfloat[MNIST \cite{lecun2010mnist}]{\label{subfig:class1}\includegraphics[width=0.3\textwidth]{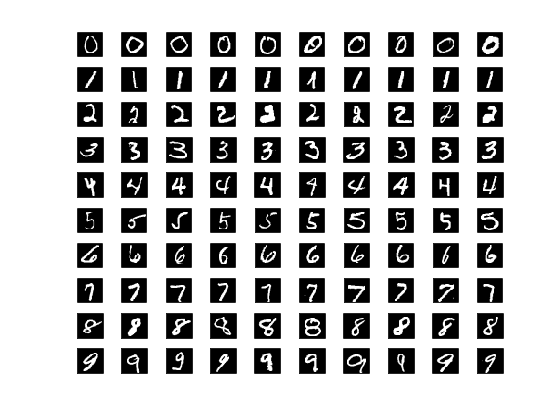}}\vspace{0.05cm}
\subfloat[Prediction]{\label{subfig:class2}\includegraphics[width=0.3\textwidth,trim = 0cm 7cm 0cm 7cm]{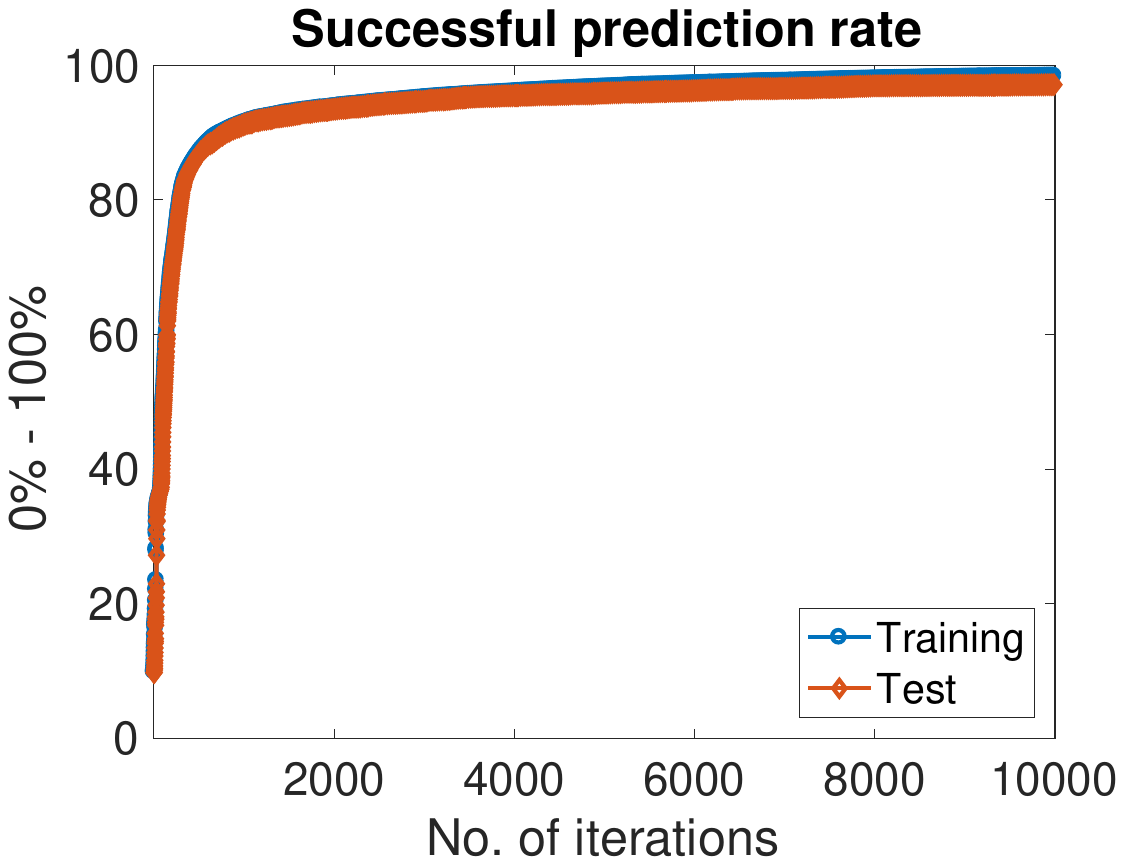}}\vspace{0.05cm}
\subfloat[Rank]{\label{subfig:class3}\includegraphics[width=0.3\textwidth,trim = 0cm 7cm 0cm 7cm]{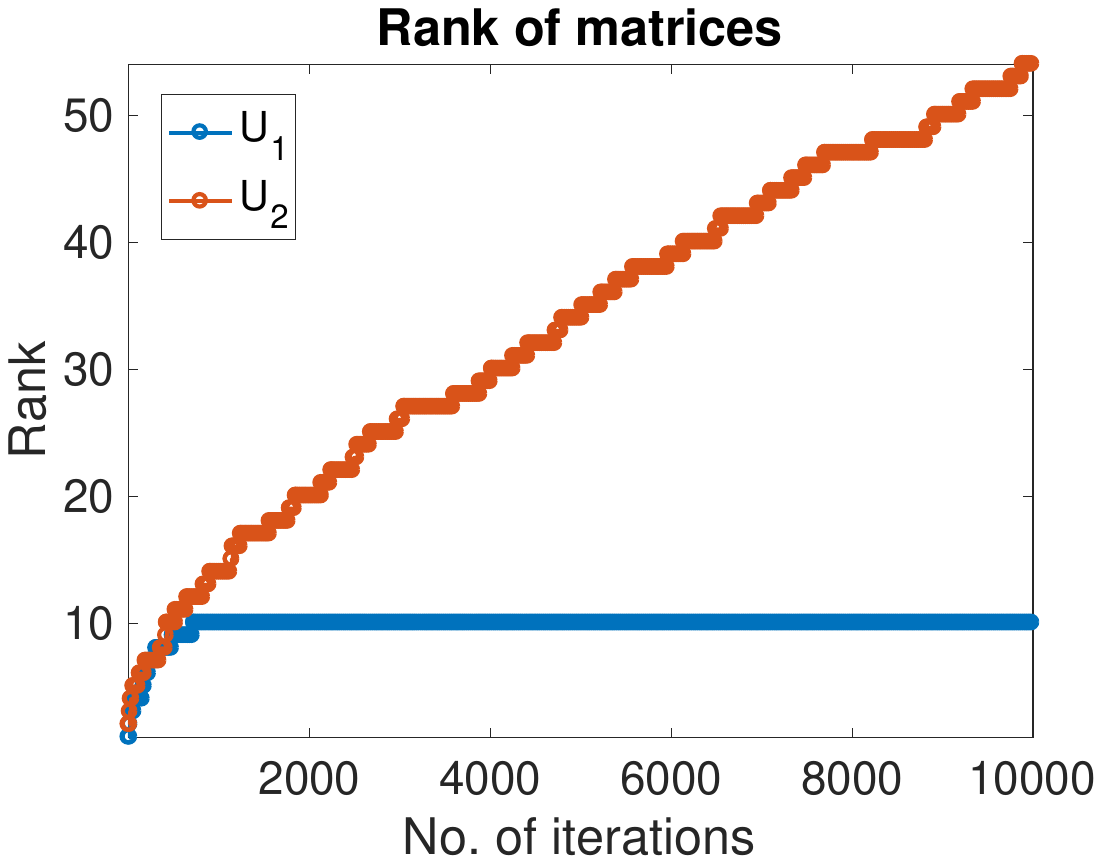}}
\end{center}
\caption{Figure \ref{subfig:class1} shows ten randomly chosen images of each digit from the MNIST training data \cite{lecun2010mnist}. Figure \ref{subfig:class2} shows the successful prediction rate of the classifier throughout the iteration both for the training and the test data. Figure \ref{subfig:class3} shows the rank of the two matrices $U_1$ and $U_2$ that are reconstructed. It becomes evident that the rank is monotonically increasing throughout the course of the iteration, allowing the model to fit only a reduced no. of effective parameters at a time.}
\end{figure}

\subsection{Classification}
We test the proposed framework for the classification of images of hand-written digits. We use the well-known MNIST dataset \cite{lecun2010mnist} as the basis for our classification. Ten example images of each class are visualised in Figure \ref{subfig:class1}. We pick 50000 images from the training dataset to create our training data matrix $D$, and use the remaining 10000 for cross validation. We model our classifier as a two-level neural network as described in Section \ref{subseq:classification}. We choose the original rectifier activation functions for the networks' architecture, in order to ensure that the composition is semi-algebraic and that the KL condition is satisfied. We overcome the non-differentiability by setting the derivatives to zero at the non-differentiable points. This is consistent with the smooth-$\max$ approximation of the rectifier for $\beta \rightarrow \infty$. We choose $E$ to be the squared Frobenius norm and set the scaling parameters to $\alpha_1 = \alpha_2 = 0.2$. The stepsize $\tau^0$ is initialised with $\tau^0 = 10^{-3}$. Subsequently, we run Algorithm \ref{alg:linbregspecial} for 10000 iterations. The prediction results of the classifier and the rank of the trained matrices are visualised in Figure \ref{fig:classification}.

\section{Conclusions \& Outlook}\label{sec:concoutlook}
We have presented a generalisation of gradient descent that allows the incorporation of non-smooth Bregman distances, and therefore can also be seen as an extension of the linearised Bregman iteration to non-convex functions. We have shown that the proposed method satisfies a sufficient decrease property and that the computed subgradients are bounded by the gap of the primal iterates. We have proven a global convergence result, where the limit is guaranteed to be a critical point of the energy if the subgradients are locally also bounded. The numerical experiments suggest that the proposed method together with early stopping can be designed to obtain solutions superior to those attained with conventional variational regularisation methods.

There are several open questions and natural directions that can be explored from here. One could extend the method to more general proximal mappings, as demonstrated in an earlier preprint. One could also study a linearised block coordinate variant of the proposed method, which would be similar in analysis to \cite{xu2013block,bolte2014proximal}. In the wake of \cite{ochs2014ipiano,pock2016inertial}, a generalisation of the proposed method could include inertial terms (or even multi-step inertial terms as in \cite{liang2016multi}), or Nesterov acceleration as in \cite{huang2013accelerated}. Both approaches seem intuitive for accelerating the method. Another direction that can be explored is the direction of non-smooth quasi-Newton extensions similar to \cite{becker2012quasi}. Motivated by applications in deep learning, one could also follow up on incremental or stochastic variants of the proposed algorithm (cf. \cite{Johnson2013,Defazio2014,Bertsekas2011a}). As we have used early stopping in our practical experiments, an interesting open question is whether the linearised Bregman iteration is a regularisation method, and if so, in what sense. This has been partially addressed in \cite{Bachmayr2009}, but under more restrictive assumptions. Following diagonal iterative regularisation approaches, an interesting open question is also if the concept of \cite{garrigos2016iterative} can be combined with the linearised Bregman iteration for non-convex problems.

\bibliographystyle{siamplain}
\bibliography{bib.bib}

\appendix
\section{Mathematical preliminaries}\label{sec:mathprelim}
We briefly summarise several concepts of convex and non-convex analysis that are of importance for the main part of this paper. Detailed informations about these concepts can be found in various textbooks, such as \cite{rockafellar1970convex,bauschke2011convex}. We frequently use functions that are proper, lower semi-continuous and convex, and therefore define the following set of functions:
\begin{align*}
\FunConvex := \left\{ J: \uc \rightarrow \R \cup \{ \infty \} \, | \, \text{$J$ is proper, lower semi-continuous and convex} \right\} \, .
\end{align*}
Here proper means that the effective domain of $J$ is not empty. The effective domain of $J$ is defined as follows.
\begin{definition}[Effective domain]
The \emph{effective domain} of a function $J:\R^n \rightarrow \R \cup \{ \infty \}$ is defined as 
\begin{align*}
\dom(J) := \{ u \in \R^n \, | \, J(u) < \infty \} \, .
\end{align*}
\end{definition}
Convex and proper functions are not necessarily differentiable, but subdifferentiable. We therefore want to recall the definition of subgradients and the subdifferential of a convex function.
\begin{definition}[Subdifferential]
Let $J \in \FunConvex$. The function $J$ is called \emph{subdifferentiable} at $u \in \uc$, if there exists an element $p \in \uc$ such that
\begin{align*}
J(v) \geq J(u) + \langle p, v - u \rangle 
\end{align*}
holds, for all $v \in \uc$. Furthermore, we call $p$ a \emph{subgradient} at position $u$. The collection of all subgradients at position $u$, i.e.
\begin{align*}
\partial J(u) := \left\{ p \in \uc \, | \, J(v) \geq J(u) + \langle p, v - u \rangle \, , \, \forall v \in \uc \right\} \, ,
\end{align*}
is called \emph{subdifferential} of $J$ at $u$.
\end{definition}
Another useful concept that we want to recall is the concept of Fenchel-, respectively convex-conjugates.
\begin{definition}[Convex conjugate]\label{def:convconj}
Let $J \in \FunConvex$. Then its \emph{convex conjugate} $J^\ast: \R^n \rightarrow \R \cup \{ \infty \}$ is defined as 
\begin{align*}
J^\ast(p) := \sup_{u \in \R^n} \left\{ \langle u, p \rangle - J(u) \right\} \, ,
\end{align*}
for all $p \in \R^n$.
\end{definition}
Amongst others, subgradients of convex conjugates satisfy the following two useful properties.
\begin{lemma}\label{lem:convconj}
Let $J \in \FunConvex$, and $J^\ast$ denote the convex conjugate of $J$. Then for all arguments $u \in \R^n$ with corresponding subgradients $p \in \partial J(u)$ we know
\begin{itemize}
\item $\langle u, p \rangle = J(u) + J^\ast(p)$,
\item $p \in \partial J(u)$ is equivalent to $u \in \partial J^\ast(p)$.
\end{itemize}
\end{lemma}

Bregman distances, introduced by Lev Bregman in 1967 (see \cite{bregman1967relaxation}), play a vital role in the definition as well as in the convergence analysis of the linearised Bregman iteration for non-convex functions. We recall its generalised variant for subdifferentiable functions \cite{kiwiel1997proximal}.

\begin{definition}[Bregman distance]
Let $J \in \FunConvex$. Then the \emph{generalised Bregman distance} for a particular subgradient $q \in \partial J(v)$ is defined as
\begin{align}
D_J^q(u, v) := J(u) - J(v) - \langle q, u - v \rangle \, , \label{eq:bregdis}
\end{align}
for $v \in \dom(J)$ and all $u \in \R^n$. 
\end{definition}


\begin{remark}\label{rem:bregdis}
Based on Lemma \ref{lem:convconj} we can rewrite \eqref{eq:bregdis} as follows:
\begin{align}
D^q_J(u, v) = J(u) + J^\ast(q) - \langle u, q \rangle \, .\label{eq:bregdis2}
\end{align}
Noticeable, the Bregman distance does not depend on $v$ anymore, and could therefore be defined as a function of $u$ and $q$ only, $D_J(u, q)$, via \eqref{eq:bregdis2} instead.
\end{remark}



Bregman distances are not symmetric in general; however, they satisfy a dual symmetry $D_J^q(u, v) = D_{J^\ast}^u(q, p)$ for arguments $u \in \uc$, $v \in \dom(J)$ and subgradients $p \in \partial J(u)$ and $q \in \partial J(v)$. Symmetry can nevertheless be achieved by simply adding two Bregman distances with interchanged arguments. The name symmetric Bregman distance goes back to \cite{burger2007error}.

\begin{definition}[Symmetric Bregman distance]
Let $J \in \FunConvex$. Then the \emph{symmetric generalised Bregman distance} $D_J^{\text{symm}}(u, v)$ is defined as
\begin{align*}
D_J^{\text{symm}}(u, v) := D_J^q(u, v) + D_J^p(v, u) = \langle p - q, u - v \rangle \, \text{,}
\end{align*}
for $u, v \in \dom(J)$ with $p \in \partial J(u)$ and $q \in \partial J(v)$.
\end{definition}


Another concept that we exploit is Lipschitz-continuity of the gradient of a function. For general operators, Lipschitz-continuity is defined as follows.

\begin{definition}[Lipschitz-continuity]\label{def:lipschitz}
An operator $F:U \subset \uc \rightarrow \R^m$ is said to be \emph{(globally) Lipschitz-continuous} if there exists a constant $L \geq 0$ such that
\begin{align}
\| F(u) - F(v) \| \leq L \| u - v \| \label{eq:lipschitz}
\end{align}
is satisfied for all $u, v \in U$. 
\end{definition}

Due to the importance of Lipschitz-continuous gradients, we define the following class of continuously differentiable functions with Lipschitz-continuous gradient:

\changed{
\begin{definition}[Smoothness]\label{def:smoothness}
A function $J : U \subset \R^n \rightarrow \R$ is called \emph{$L$-smooth} if it is differentiable and its gradient $\nabla J : U \rightarrow \uc$ is Lipschitz-continuous with Lipschitz constant $L$. The set of all $L$-smooth functions is therefore denoted by $\FunSmooth{L}$ with
\begin{align*}
\FunSmooth{L} := \left\{ J:U \rightarrow \R \left| \begin{array}{c} \text{$J$ is continuously differentiable}\\ \text{$\nabla J$ is $L$-Lipschitz-continuous} \end{array} \right. \right\} \, .
\end{align*}
\end{definition}}

\noindent Note that it is a well-known fact that $L$-smooth functions satisfy the Lipschitz estimate
\begin{align}
 J(u) \leq J(v) + \langle \nabla J(v), u - v \rangle + \frac{L}{2} \| u - v \|_2^2 \, , \label{eq:lipschitzest3}
\end{align}
\changed{for all $u, v \in U$. Note that if $U = \R^n$ then $J$ is already globally $L$-smooth and this estimate is true for all arguments $u, v \in \R^n$.}

\noindent In the following we recall the definition of the proximal mapping.
\begin{definition}[Proximal mapping \cite{moreau1962decomposition,moreau1965proximite}]
We define the \emph{proximal mapping} as the operator $(I + \partial J)^{-1}:\uc \rightarrow \dom(J)$ with
\begin{align*}
(I + \partial J)^{-1}( f ) := \argmin_{u \in \dom(J)} \left\{ \frac{1}{2}\| u - f \|^2 + J(u) \right\} \, ,
\end{align*}
for all arguments $f \in \uc$.
\end{definition}

To conclude this section, we want to recall the Kurdyka-\L ojasiewicz (KL) property \cite{lojasiewicz1963propriete,kurdyka1998gradients}. For the definition of the KL property we need to define a distance between sub-sets and elements of $\R^n$ first.
\begin{definition}\label{def:distance}
Let $\Omega \subset \R^n$ and $u \in \R^n$. We define the distance from $\Omega$ to $u$ as
\begin{align*}
\dist(u, \Omega) := \begin{cases} \inf \{ \| v - u \| \, | \, v \in \Omega \} & \Omega \neq \emptyset\\ \infty & \Omega = \emptyset \end{cases} \, .
\end{align*}
\end{definition}
The definition of the KL property based on the distance measure defined in Definition \ref{def:distance} reads as follows.

\begin{definition}[Kurdyka-\L ojasierwicz property]\label{def:kl}
A function $J$ is said to have the \emph{Kurdyka-\L ojasierwicz (KL) property} at $\overline{u} \in \text{dom}(\partial J) := \{ u \in \R^n \, | \, \partial J(u) \neq \emptyset \}$ if there exists a constant $\eta \in (0, \infty]$, a neighbourhood $\Theta$ of $\overline{u}$ and a function $\varphi:[0, \eta ) \rightarrow \R_{> 0}$, which is a concave function that is continuous at 0 and satisfies $\varphi(0) = 0$, $\varphi \in C^1((0, \eta))$ and $\varphi^\prime(s) > 0$ for all $s \in (0, \eta)$, such that for all $u \in \Theta \cap \{ u \in \R^n \, | \, J(\overline{u}) < J(u) < J(\overline{u}) + \eta \}$ the inequality
\begin{align}
\varphi^\prime\left(J(u) - J(\overline{u})\right)\dist(0, \partial J(u)) \geq 1 \tag{KL} \label{eq:modifiedkl}
\end{align}
holds.

\noindent If $J$ satisfies the KL property at each point of $\dom(\partial J)$, $J$ is called a \emph{KL function}.
\end{definition}
We conclude the appendix by recalling one important result from \cite{bolte2014proximal} that is necessary for successfully carrying out the convergence proof in the main part of the paper.
\begin{lemma}[{Uniformised KL property \cite[Lemma 6]{bolte2014proximal}}]\label{lem:unifiedkl}
Let $\Omega$ be a compact set, and suppose that $J$ is a function that is constant on $\Omega$ and that satisfies \eqref{eq:modifiedkl} at each point in $\Omega$. Then there exist $\varepsilon > 0$, $\eta > 0$ and $\varphi \in C^1( (0, \eta) )$ that satisfy the same conditions as in Definition \ref{def:kl}, such that for all $\overline{u} \in \Omega$ and all $u$ in
\begin{align}
\left\{ u \in \R^n \, \left| \, \dist(u, \Omega) < \varepsilon \right. \right\} \cap \left\{ u \in \R^n \, \left| \, J(\overline{u}) < J(u) < J(\overline{u}) + \eta \right. \right\} \label{eq:klset}
\end{align}
condition \eqref{eq:modifiedkl} is satisfied.
\end{lemma}

\end{document}